# MULTIPLE DECORRELATION AND RATE OF CONVERGENCE IN MULTIDIMENSIONAL LIMIT THEOREMS FOR THE PROKHOROV METRIC

BY FRANCOISE PÈNE

*University of Brest*

The motivation of this work is the study of the error term $e_t^\varepsilon(x,\omega)$ in the averaging method for differential equations perturbed by a dynamical system. Results of convergence in distribution for $(\frac{e_t^\varepsilon(x,\cdot)}{\sqrt{\varepsilon}})_{\varepsilon>0}$ have been established in Khas'minskii [*Theory Probab. Appl.* **11** (1966) 211–228], Kifer [*Ergodic Theory Dynamical Systems* **15** (1995) 1143–1172] and Pène [*ESAIM Probab. Statist.* **6** (2002) 33–88]. We are interested here in the question of the rate of convergence in distribution of the family of random variables $(\frac{e_t^\varepsilon(x,\cdot)}{\sqrt{\varepsilon}})_{\varepsilon>0}$ when $\varepsilon$ goes to 0 ($t>0$ and $x \in \mathbf{R}^d$ being fixed). We will make an assumption of multiple decorrelation property (satisfied in several situations). We start by establishing a simpler result: the rate of convergence in the central limit theorem for regular multidimensional functions. In this context, we prove a result of convergence in distribution with rate of convergence in $O(n^{-1/2+\alpha})$ for all $\alpha > 0$ (for the Prokhorov metric). This result can be seen as an extension of the main result of Pène [*Comm. Math. Phys.* **225** (2002) 91–119] to the case of $d$-dimensional functions. In a second time, we use the same method to establish a result of convergence in distribution for $(\frac{e_t^\varepsilon(x,\cdot)}{\sqrt{\varepsilon}})_{\varepsilon>0}$ with rate of convergence in $O(\varepsilon^{1/2-\alpha})$ (for the Prokhorov metric). We close this paper with a discussion (in the Appendix) about the behavior of the quantity $\|\sup_{0 \le t \le T_0} |e_t^\varepsilon(x,\cdot)|_\infty\|_{L^p}$ under less stringent hypotheses.

**1. Introduction.** We are interested in the asymptotic behavior of random variables sequences defined by a probability dynamical system. Let us consider a (discrete-time) probability dynamical system $(\Omega, \mathcal{F}, \nu, T)$ [where $(\Omega, \mathcal{F}, \nu)$ is a probability space endowed with a $\nu$-preserving transformation $T:\Omega \to \Omega$].

Let a function $f$ defined on $\Omega$ with values in $\mathbf{R}^d$ be given. We can study the stochastic properties of the sequence of random variables $(f \circ T^n)_{n \geq 0}$ defined on $(\Omega, \mathcal{F}, \nu)$. If $(\Omega, \mathcal{F}, \nu, T)$ is ergodic, Birkhoff's ergodic theorem [6] gives a strong law of large numbers for $(f \circ T^n)_{n \geq 0}$ when the function $f$ is $\nu$-integrable. Furthermore, central limit theorems (CLTs) have been established for $(f \circ T^n)_{n \geq 0}$ in various situations (see [9, 31, 36, 38], etc.). Results of speed of convergence in the CLT for $(f \circ T^n)_{n \geq 0}$ have been established in the one-dimensional case (i.e., when $f$ is a real-valued function) in [17], [19] and [32], for example. Here, we are interested in the speed of convergence in the central limit theorem for multidimensional random variables $(f \circ T^n)_n$ (i.e., when $d \geq 2$). We estimate the speed in the sense of the Prokhorov metric. When $(f \circ T^n)_n$ is a sequence of independent random variables, Yurinskii established a speed of convergence in $\frac{1}{\sqrt{n}}$ in the sense of the Prokhorov metric (cf. [39]). Let us point out the fact that such an estimate gives directly an estimate in $\frac{1}{\sqrt{n}}$ for the speed of convergence of the expectation of any bounded Lipschitz continuous function. In Section 2 of the present paper we establish a speed of convergence in $O(n^{-1/2+\alpha})$ (for all $\alpha > 0$) for the multidimensional CLT for $(f \circ T^n)_{n \geq 0}$ when $f$ is a regular function (Theorem 2.2). This result holds under a hypothesis of multiple decorrelation (with exponential rate) for regular functions. This hypothesis is satisfied in different hyperbolic situations (systems studied in [38], billiard transformation studied in [37], mostly contracting diffeomorphisms studied in [8]).

Our proof is based on the method developed by Jan to establish Theorem 7 of [19] (it uses characteristic functions) and on a result due to Yurinskii [39] which plays here a similar role to the one played by the more classical Esseen lemma [12] in the proof of Theorem 7 of [19]. (Let us mention the work of Jan who estimated, in a slightly different context, the speed of convergence in the multidimensional central limit theorem in the sense of the uniform convergence of the distribution functions and then extended Rio's result of [32]; cf. Theorem 9 of [19].)

In Section 3, a result of speed of convergence in terms of the Prokhorov metric is established in a more sophisticated context. We study the averaging method for differential equations perturbed by the probability dynamical system $(\Omega, \mathcal{F}, \nu, T)$. This problem has been studied in particular [11, 20, 21, 25, 26]. For a general reference about this method, we refer to Chapter 4 of [1] and to Chapter 7 of [14] (see also Chapter 5 of [2]). The problem is the following one. Let a function $F : \mathbf{R}^d \times \Omega \to \mathbf{R}^d$ smooth enough (measurable, uniformly bounded and uniformly Lipschitz in the first parameter) be given. For any $\varepsilon > 0$ and any $(x, \omega) \in \mathbf{R}^d \times \Omega$, we consider the continuous solution $(x_t^\varepsilon(x, \omega))_t$ of the following differential equation (with initial condition):

$$\forall t \in \mathbf{R}_+ \setminus \mathbf{N}, \qquad \frac{dx_t^\varepsilon}{dt}(x, \omega) = F(x_t^\varepsilon(x, \omega), T^{\lfloor t/\varepsilon \rfloor}(\omega)) \quad \text{and} \quad x_0^\varepsilon(x, \omega) = x.$$

Let us write $(w_t(x))_t$ the solution of the differential equation (with initial condition) obtained from the previous one by averaging:

$$\forall t \in \mathbf{R}_+ \qquad \frac{dw_t}{dt}(x) = \int_\Omega F(w_t(x), \omega')\, d\nu(\omega') \quad \text{and} \quad w_0(x) = x.$$

We are interested in the study of the asymptotic behavior (when $\varepsilon$ goes to 0) of the error term $(e_t^\varepsilon(x, \omega))_t$ defined by

$$e_t^\varepsilon(x, \omega) := x_t^\varepsilon(x, \omega) - w_t(x).$$

Results of convergence in distribution for the family of processes $((\frac{e_t^\varepsilon(x,\cdot)}{\sqrt{\varepsilon}})_{t \in [0;T]})_{\varepsilon > 0}$ have been established in [21] and in [28] (see Theorem 2.1.3 in [28]). Here, we establish a result of speed of convergence in distribution for the family of random variables $(\frac{e_s^\varepsilon(x,\cdot)}{\sqrt{\varepsilon}})_{\varepsilon > 0}$, $s$ and $x$ being fixed (Theorem 3.4). The speed is estimated in the sense of the Prokhorov metric. The proof of this result is based on the ideas of the proof of Theorem 2.2.

In the Appendix, we complete our study with estimates of the following form:

$$\sup_{x \in \mathbf{R}^d} \left\| \sup_{t \in [0;T_0]} |e_t^\varepsilon(x,\cdot)|_\infty \right\|_{L^p} = O(\sqrt{\varepsilon}),$$

for any real number $T_0 > 0$ and for some real number $p \geq 1$. With these results, we improve a result of [11] in two particular cases: for a differential equation perturbed by the billiard flow studied in [29, 37] and in the case of a differential equation perturbed by a diagonal flow on a compact quotient of $SL(d, \mathbf{R})$ (see Section A.1.4).

1.1. *Context.* Let us specify the context we consider here. Let us consider a probability dynamical system $(\Omega, \mathcal{F}, \nu, T)$. Let us suppose that the space $\Omega$ is endowed with a metric $d$ and that $\mathcal{F}$ is the associated Borel $\sigma$-algebra. We denote by $\mathbf{E}_\nu[\cdot]$ the expectation relative to the measure $\nu$:

$$\mathbf{E}_\nu[f] := \int_\Omega f\, d\nu.$$

For all complex-valued square integrable functions $f, g$, we denote by $\mathrm{Cov}_\nu(f, g)$ the covariance of the functions $f$ and $g$ with respect to the measure $\nu$:

$$\mathrm{Cov}_\nu(f, g) = \mathbf{E}_\nu[fg] - \mathbf{E}_\nu[f]\mathbf{E}_\nu[g].$$

Let a real number $\eta \in ]0; 1]$ be fixed. For any uniformly bounded and $\eta$-Hölder continuous function $f : \Omega \to \mathbf{C}$, we define $\|f\|_\infty := \sup_{x \in \Omega} |f(x)|$ and we denote by $C_f^{(\eta)}$ the Hölder coefficient of order $\eta$ of $f$:

$$C_f^{(\eta)} := \sup_{x \neq y} \frac{|f(x) - f(y)|}{d(x,y)^\eta}.$$

We write $\mathcal{H}_\eta$ the set of complex-valued uniformly bounded $\eta$-Hölder continuous functions defined on $\Omega$.

For any real number $r \geq 1$, we introduce the multiple decorrelation Property $(\mathcal{P}_r)$ as follows:

PROPERTY $(\mathcal{P}_r)$. *There exist a polynomial function $P_r$ with real nonnegative coefficients and a real number $\delta_r \in ]0;1[$ such that, for all integers $m$ and $m'$, for all bounded $\eta$-Hölder continuous functions $f_1, \ldots, f_{m+m'} : \Omega \to \mathbf{C}$, for all increasing finite sequences of nonnegative integers $(k_1, \ldots, k_m)$ and $(l_1, \ldots, l_{m'})$ and for all nonnegative integer $n$, we have*

$$
(1) \quad \left| \operatorname{Cov}_\nu \left( \prod_{i=1}^m f_i \circ T^{k_i}, \prod_{j=1}^{m'} f_{m+j} \circ T^{n+l_j} \right) \right|
$$
$$
\leq \left( \prod_{i=1}^{m+m'} \|f_i\|_\infty + \sum_{i=1}^{m+m'} C^{(\eta)}_{f_i} \prod_{j \neq i} \|f_j\|_\infty \right) P_r(l_{m'}) \delta_r^{n-rk_m}.
$$

Such results of decorrelation have been studied in [22] for Anosov diffeomorphisms. Let us make some commentaries about this property. Let us notice that Theorems 2.2 and 3.4 are still true if we replace, in Property $(\mathcal{P}_r)$, $\delta_r^{n-rk_m}$ by $h_r(n - rk_m)$, where $(h_r(n))_{n \geq 0}$ decreases rapidly (more precisely, if $\lim_{n \to +\infty} n^\beta h_r(n) = 0$ for every real number $\beta > 0$).

Property $(\mathcal{P}_r)$ is satisfied for any $r > 1$ in the case of a billiard transformation studied in [37] (cf. Corollary B.2. of [29]). This result can be proved in the same way for any dynamical system to which Young's method of [38] can be applied. Examples of dynamical systems satisfying this property are given in [23] where a similar property is proved. In particular, this property is satisfied for ergodic algebraic automorphisms of the torus (this can be proved by rewriting the proof of Theorem 4.1.2 of [28]) and for diagonal transformation on a compact quotient of $SL(d, \mathbf{R})$ (see [23]) and for the dynamical systems studied by Dolgopyat in [9].

1.2. *Prokhorov metric, definition and first results.* We endow $\mathbf{R}^d$ with the supremum norm $|\cdot|_\infty$ defined by $|(x_1, \ldots, x_d)|_\infty := \max_{i=1, \ldots, d} |x_i|$. For real-valued random variables, we estimate the speed of convergence in distribution in terms of uniform convergence of distribution functions. In the $d$-dimensional case, a natural metric between two probability measures on $\mathbf{R}^d$ is the Prokhorov metric (cf. [10], e.g.). Let us recall now its definition and some of its properties.

DEFINITION 1.1 (Prokhorov metric). Let $P$ and $Q$ be two probability measures on $\mathbf{R}^d$. The Prokhorov metric $\Pi(P, Q)$ between $P$ and $Q$ is the

following quantity:

$$\Pi(P,Q) := \inf\left\{\varepsilon > 0 : \sup_{B \in \mathcal{B}}(P(B) - Q(B^\varepsilon)) \leq \varepsilon\right\},$$

where $\mathcal{B}$ is the Borel $\sigma$-algebra on $\mathbf{R}^d$ and where we denote by $B^\varepsilon$ the $\varepsilon$-open neighborhood of $B$.

Let us recall the link between the Prokhorov metric for the probability measures on $\mathbf{R}^d$ and the Ky Fan metric for the $\mathbf{R}^d$-valued random variables defined on the same probability space.

DEFINITION 1.2 (Ky Fan metric). Let $X$ and $Y$ be two $\mathbf{R}^d$-valued random variables defined on the same probability space $(E, \mathcal{T}, \mathbf{P})$. The Ky Fan metric (associated to $|\cdot|_\infty$) between $X$ and $Y$ is given by

$$\mathcal{K}(X,Y) := \inf\{\varepsilon > 0 : \mathbf{P}(|X - Y|_\infty > \varepsilon) < \varepsilon\}.$$

PROPOSITION 1.3. *Let $P$ and $Q$ be two Borel probability measures on $\mathbf{R}^d$. The Prokhorov metric $\Pi(P,Q)$ between $P$ and $Q$ is the infimum of the Ky Fan metric between $X$ and $Y$, where $(X, Y)$ describes the set of couples of random variables defined on the same probability space such that the distribution of $X$ is $P$ and such that the distribution of $Y$ is $Q$.*

Another classical metric between probability measures on $\mathbf{R}^d$ is the BL metric (BL for *bounded Lipschitz*) defined as follows:

DEFINITION 1.4. Let $P$ and $Q$ be two probability measures on $\mathbf{R}^d$. The BL metric between $P$ and $Q$ is the following quantity:

$$BL(P,Q) := \sup\left\{\frac{\mathbf{E}_P[\phi] - \mathbf{E}_Q[\phi]}{\|\phi\|_\infty + L_\phi} \bigg| \phi : \mathbf{R}^d \to \mathbf{R},\ \|\phi\|_\infty + L_\phi < +\infty\right\},$$

where we denote $\|\phi\|_\infty = \sup_{x \in \mathbf{R}^d} |\phi(x)|$ and $L_\phi := \sup_{x \neq y} \frac{|\phi(x) - \phi(y)|}{|x-y|_\infty}$.

These two metrics are metrics for the weak convergence for probability measures (which corresponds to the convergence in distribution for random variables). Moreover, we have the following (cf., e.g., [24], Proposition 1.2 and [10], Problem 11.3.5):

PROPOSITION 1.5 (Equivalence of these metrics). *Let $P$ and $Q$ be two Borel probability measures on $\mathbf{R}^d$. We have*

$$\tfrac{1}{3}BL(P,Q) \leq \Pi(P,Q) \leq (\tfrac{3}{2}BL(P,Q))^{1/2}.$$

In the following, we will essentially be interested in questions of speed of convergence in terms of Prokhorov metric. But, we will also talk about BL metric.

1.3. *Notation.* Let $A$ and $B$ be any vectors in $\mathbf{R}^d$. Let us denote by $^TA$ the line vector, transposed to $A$. Let us denote by $A \otimes B$ the square $d$-dimensional matrix given by $A \otimes B := A \cdot {}^TB$ and we write $A^{\otimes 2} := A \otimes A$.

Let a probability space $(\Omega, \mathcal{F}, \nu)$ and a real number $p \geq 1$ be given. We denote by $L^p(\Omega, \mathbf{R}^d)$ the set of measurable functions $f : \Omega \to \mathbf{R}^d$ such that $\int_\Omega |f|_\infty^p \, d\nu < +\infty$. For any $y$ in $L^p(\Omega, \mathbf{R}^d)$, we denote $\|f\|_{L^p} = (\int_\Omega |f|_\infty^p \, d\nu)^{1/p}$.

For any probability space $(\Omega, \mathcal{F}, \nu)$, any measurable space $(E, \mathcal{T})$ and any random variable $X : \Omega \to E$, we denote by $\nu_*(X)$ the image measure of $\nu$ by $X$, that is, the probability measure defined on $(E, \mathcal{T})$ by $\nu_*(X)(A) = \nu(X^{-1}(A))$ for any $A \in \mathcal{T}$.

## 2. Ordinary central limit theorem.

2.1. *Introduction and result in the i.i.d. case.* We are interested here in the question of the rate of convergence in the central limit theorem, that is, the question of the rate of convergence in distribution for sequences of random variables of the form $(\frac{1}{\sqrt{n}} \sum_{k=0}^{n-1} X_k)_{n \geq 1}$ to a normal random variable. For any $A \in \mathbf{R}^d$ and any $d \times d$ nonnegative symmetric matrix $C$, we denote by $\mathcal{N}(A, C)$ the normal distribution with mean $A$ and with covariance matrix $C$ (cf. [13], III-6, for the notion of normal distributions).

For independent multidimensional variables, results of speed of convergence have been established by many authors under moment hypotheses. Let us mention the works of Bergstrm [3], Sazanov [34], Ranga Rao [30] and Bhattacharya [4] (for uniform estimates) and of Rotar [33] (for a nonuniform estimate). Let us give the following result coming from [39]. The proof of this result given by Yurinskii is based on a result linking Prokhorov metric with characteristic functions (cf. Proposition 2.6).

THEOREM 2.1. *Let $(X_k)_{k \geq 0}$ be a sequence of $\mathbf{R}^d$-random variables defined on a probability space $(\Omega, \mathcal{F}, \mathbf{P})$. If these random variables are independent and identically distributed, $\mathbf{P}$-centered and admitting moments of the third order, then the sequence of random variables $(\frac{1}{\sqrt{n}} \sum_{k=0}^{n-1} X_k)_{n \geq 1}$ converges in distribution to a random variable with (eventually degenerate) normal distribution $\mathcal{N}(0, \mathbf{E}[X_1^{\otimes 2}])$ and we have*

$$\Pi\left(\mathbf{P}_*\left(\frac{1}{\sqrt{n}} \sum_{k=0}^{n-1} X_k\right), \mathcal{N}(0, \mathbf{E}[X_1^{\otimes 2}])\right) = O\left(\frac{1}{\sqrt{n}}\right).$$

Moreover, this speed of convergence is optimal under these hypotheses [there exists such a sequence of random variables $(X_k)_k$ for which the speed is exactly in $\frac{1}{\sqrt{n}}$].

Here we will consider random variables $X_k$ which are maybe not independent but are stationary. More precisely, we will suppose that the random variables $X_k$ are given by $X_k = f \circ T^k$ with $(\Omega, \mathcal{F}, \nu, T)$ as described before and with $f : \Omega \to \mathbf{R}^d$ any uniformly bounded $\eta$-Hölder continuous function.

2.2. *A rate of convergence in the central limit theorem.* For any function $f : \Omega \to \mathbf{R}^d$ and any integer $n \geq 1$, we define

$$S_n(f) := \sum_{k=0}^{n-1} f \circ T^k.$$

First of all, let us notice that, under hypothesis $(\mathcal{P}_r)$, the following limit exists for any $\nu$-centered, bounded $\eta$-Hölder continous function $f : \Omega \to \mathbf{R}^d$:

$$D(f) := \lim_{n \to +\infty} \left( \mathbf{E}_\nu \left[ \left( \frac{S_n(f)}{\sqrt{n}} \right)^{\otimes 2} \right] \right),$$

and that we have

(2) $\qquad D(f) = \mathbf{E}_\nu[f \otimes f] + \sum_{k \geq 1} (\mathbf{E}_\nu[f \otimes f \circ T^k] + \mathbf{E}_\nu[f \circ T^k \otimes f]).$

THEOREM 2.2. *We suppose that there exists some $r \geq 1$ for which Property $(\mathcal{P}_r)$ is satisfied. Let $f : \Omega \to \mathbf{R}^d$ be a $\nu$-centered, bounded $\eta$-Hölder continous function. If the matrix $D(f)$ is nondegenerate, then the sequence of random variables $(\frac{S_n(f)}{\sqrt{n}})_{n \geq 0}$ converges in distribution to a $d$-dimensional random variable with normal distribution $\mathcal{N}(0, D(f))$ and we have*

(3) $\quad \forall \alpha > 0, \qquad \Pi_n(f) := \Pi\left( \nu_* \left( \frac{1}{\sqrt{n}} S_n(f) \right), \mathcal{N}(0, D(f)) \right) = O(n^{-1/2+\alpha}).$

Let us make some comments about the case in which the asymptotic covariance $D(f)$ is degenerate. By a classical argument (cf., e.g., Lemma 2.2 of [7]), we have the following result:

PROPOSITION 2.3. *Let us suppose that there exists a real number $r \geq 1$ for which Property $(\mathcal{P}_r)$ is satisfied. If $g : \Omega \to \mathbf{R}$ is a $\nu$-centered, bounded $\eta$-Hölder continuous function such that $D(g) = 0$, then $g$ is a coboundary in $L^2$, that is, there exists a $\nu$-centered square integrable function $h : \Omega \to \mathbf{R}$ such that we have $g = h - h \circ T$ almost surely.*

If $f : \Omega \to \mathbf{R}^d$ is a $\nu$-centered, bounded $\eta$-Hölder continuous function, then there exists an orthogonal matrix $A \in O_d(\mathbf{R})$ such that the matrix $D(A \cdot f) = A \cdot D(f) \cdot {}^T A$ is diagonal with diagonal terms $\alpha_1 \geq \alpha_2 \geq \cdots \geq \alpha_d$. Let us suppose now that the matrix $D(f)$ is degenerate. It is natural to ask if, in that case, estimate (3) is still true. Because of the equivalence of norms in finite dimension, estimate (3) will be true for $f$ if and only if it is true for the function $g$ defined by $g := A \cdot f$. Let $r$ be the rank of the matrix $D(f)$ and $g_1, \ldots, g_d$ be the coordinate functions of $g$. Coefficients $\alpha_1, \ldots, \alpha_r$ are nonnull positive and coefficients $\alpha_{r+1}, \ldots, \alpha_d$ are null. We can therefore apply Theorem 2.2 to the function $(g_1, \ldots, g_r) : \Omega \to \mathbf{R}^r$ and, consequently, to the function $G = (g_1, \ldots, g_r, 0, \ldots, 0) : \Omega \to \mathbf{R}^d$. Hence, we have $g = G + H$ with $D(H) = 0$. Then, according to the previous proposition applied to the coordinate functions of $H$, there exists a $\nu$-centered square integrable function $h : \Omega \to \mathbf{R}^d$ such that we have $g = G + h - h \circ T$ almost surely. Therefore, for any integer $n \geq 1$, we have

$$\frac{1}{\sqrt{n}} S_n(g) = \frac{S_n(G) + B_n}{\sqrt{n}},$$

where $(B_n := h - h \circ T^n)_n$ is a sequence of random variables bounded in $L^2$ and with

$$\forall \alpha > 0, \qquad \Pi\left(\nu_*\left(\frac{1}{\sqrt{n}} S_n(G)\right), \mathcal{N}(0, D(g))\right) = O(n^{-1/2 + \alpha}).$$

REMARK 2.4. If the sequence of random variables $(B_n)_n$ is bounded in $L^p$ (for some $p \geq 1$), then we have (according to Markov's inequality)

$$\mathcal{K}\left(\frac{S_n(G) + B_n}{\sqrt{n}}, \frac{S_n(G)}{\sqrt{n}}\right) \leq \frac{\sup_m \|B_m\|_{L^p}^{p/(p+1)}}{n^{p/(2(p+1))}},$$

and therefore, according to Theorem 2.2,

$$\Pi\left(\nu_*\left(\frac{1}{\sqrt{n}} S_n(f)\right), \mathcal{N}(0, D(f))\right) = O(n^{-p/(2(p+1))}).$$

If $(B_n)_n$ is bounded in $L^p$ for all real number $p \geq 1$, then we have

$$\forall \alpha > 0, \qquad \mathcal{K}\left(\frac{S_n(G) + B_n}{\sqrt{n}}, \frac{S_n(G)}{\sqrt{n}}\right) = O(n^{-1/2 + \alpha}),$$

and therefore, according to Theorem 2.2,

$$\forall \alpha > 0, \qquad \Pi\left(\nu_*\left(\frac{1}{\sqrt{n}} S_n(f)\right), \mathcal{N}(0, D(f))\right) = O(n^{-1/2 + \alpha}).$$

If $(B_n)_n$ is bounded in $L^1$, then for any bounded Lipschitz continuous function $\phi : \mathbf{R}^d \to \mathbf{R}$, we have

$$\left|\mathbf{E}_\nu\left[\phi\left(\frac{S_n(G) + B_n}{\sqrt{n}}\right)\right] - \mathbf{E}_\nu\left[\phi\left(\frac{S_n(G)}{\sqrt{n}}\right)\right]\right| \leq L_\phi \frac{\sup_m \|B_m\|_{L^1}}{\sqrt{n}},$$

and therefore, according to Theorem 2.2,

$$\forall \alpha > 0, \qquad BL\left(\nu_*\left(\frac{1}{\sqrt{n}} S_n(f)\right), \mathcal{N}(0, D(f))\right) = O(n^{-1/2 + \alpha}).$$

CONSEQUENCE 2.5 (Case eventually degenerate). Let us suppose that there exists a real number $r \geq 1$ such that Property $(\mathcal{P}_r)$ is satisfied. Let $f : \Omega \to \mathbf{R}^d$ be a $\nu$-centered, bounded $\eta$-Hölder continuous function. Then, we have

$$\Pi\left(\nu_*\left(\frac{1}{\sqrt{n}}S_n(f)\right), \mathcal{N}(0, D(f))\right) = O(n^{-1/3})$$

and

$$\forall \alpha > 0, \qquad BL\left(\nu_*\left(\frac{1}{\sqrt{n}}S_n(f)\right), \mathcal{N}(0, D(f))\right) = O(n^{-1/2+\alpha}).$$

2.3. *Proof.* In this section, we prove Theorem 2.2. This proof is inspired by [29]. It uses a method developped by Jan in another context (cf. [18, 19]). In order to estimate the speed of convergence in terms of the Prokhorov metric, we will use the following result:

PROPOSITION 2.6 ([39]). *Let $Q$ be a normal (nondegenerate) $d$-dimensional distribution. There exist two real numbers $c_0 > 0$ and $\Gamma > 0$ such that, for any real number $U > 0$ and for any Borel probability measure $P$ on $\mathbf{R}^d$ admitting moments of order $\lfloor d/2 \rfloor + 1$, we have*

$$\Pi(P, Q)$$
$$\leq c_0 \left[ \frac{1+\Gamma}{U} \right.$$
$$\left. + \left( \int_{|t|_\infty < U} \sum_{k=0}^{\lfloor d/2 \rfloor + 1} \sum_{\{j_1,\ldots,j_k\} \in \{1,\ldots,d\}^k} \left| \frac{\partial^k}{\partial t_{j_1} \cdots \partial t_{j_k}}(\varphi_P - \varphi_Q)(t) \right|^2 dt \right)^{1/2} \right],$$

*where we denote by $\varphi_P$ and $\varphi_Q$ the characteristic functions of the distributions $P$ and $Q$, respectively:*

$$\forall t \in \mathbf{R}^d \qquad \varphi_P(t) = \mathbf{E}_P[e^{i\langle t, \cdot \rangle}] \quad and \quad \varphi_Q(t) = \mathbf{E}_Q[e^{i\langle t, \cdot \rangle}],$$

*with $\langle \cdot, \cdot \rangle$ the usual scalar product on $\mathbf{R}^d$.*

This result links the speed of convergence in terms of the Prokhorov metric with a problem of estimation of the characteristic functions. It will play the same role in our proofs as the one played by the Esseen lemma in the proof of unidimensional central limit theorems established in [19, 29].

Let us suppose that the hypotheses of Theorem 2.2 are satisfied. Let us consider a real number $r_0 \geq 1$ such that Property $(\mathcal{P}_{r_0})$ is satisfied. Let us suppose that the matrix $D(f)$ is nondegenerate. For any $t \in \mathbf{R}^d$ and any integer $n \geq 1$, we define

$$h_n(f, t) := \mathbf{E}_\nu\left[\exp\left\{\frac{i\langle t, S_n(f)\rangle}{\sqrt{n}}\right\}\right] - \exp\left\{-\frac{\langle t, D(f)t \rangle}{2}\right\}.$$

The remainder of this section is essentially devoted to the proof of the following result. Let a real number $\alpha \in ]0; \frac{1}{2}[$ be given.

PROPOSITION 2.7. *For any integer $p \geq 0$, there exist a real number $L_p = L_{p,\alpha} > 0$ and a nonnegative functions sequence $(a_{n,p,\alpha})_{n \geq 1}$ satisfying the following:*

$$\left( \int_{|t|_\infty \leq n^{1/2-\alpha}} (a_{n,p,\alpha}(t))^2 \, dt \right)^{1/2} = O_{n \to +\infty}\left( \frac{1}{n^{1/2-\alpha}} \right),$$

*and such that, for any integer $n \geq 1$ and any $t \in \mathbf{R}^d$ satisfying $|t|_\infty \leq n^{1/2-\alpha}$, we have*

$$(4) \quad \sum_{k=0}^{\lfloor d/2 \rfloor + 1} \sum_{1 \leq j_1,\ldots,j_k \leq d} \left| \frac{\partial^k}{\partial t_{j_1} \cdots \partial t_{j_k}} h_n(f,t) \right| \leq L_p \frac{|t|_\infty^p}{n^{p(1/2-\alpha)}} + a_{n,p,\alpha}(t).$$

PROOF OF PROPOSITION 2.7. Let us prove inductively on $p$ that the following Property $(\mathcal{H}_p)$ is satisfied for all integer $p \geq 0$. $\square$

PROPERTY $(\mathcal{H}_p)$. *For any real number $\beta > 0$, there exist a real number $L_{p,\alpha,\beta} > 0$ and a sequence $(a_{n,p,\alpha,\beta}(\cdot))_n$ of nonnegative uniformly bounded functions $(a_{n,p,\alpha,\beta})_{n \geq 1}$ satisfying*

$$\limsup_{n \to +\infty} n^{1/2-\alpha} \left( \int_{|t|_\infty \leq n^{1/2-\alpha}} (1 + |t|_\infty^\beta)(a_{n,p,\alpha,\beta}(t))^2 \, dt \right)^{1/2} < +\infty$$

*and such that, for any integer $n \geq 1$ and any $t \in \mathbf{R}^d$ satisfying $|t|_\infty \leq n^{1/2-\alpha}$, we have*

$$\sum_{k=0}^{\lfloor d/2 \rfloor + 1} \sum_{1 \leq j_1,\ldots,j_k \leq d} \left| \frac{\partial^k}{\partial t_{j_1} \cdots \partial t_{j_k}} h_n(f,t) \right| \leq L_{p,\alpha,\beta} \frac{1 + |t|_\infty^p}{n^{p(1/2-\alpha)}} + a_{n,p,\alpha,\beta}(t).$$

Let us first notice that, under Property $(\mathcal{P}_r)$, for any bounded Hölder continuous function $f : \Omega \to \mathbf{R}^d$, the sequence of random variables $(\frac{S_n(f)}{\sqrt{n}})_{n \geq 1}$ is uniformly bounded in $L^p$ for any real number $p \geq 1$ (see Lemma 2.3.4 of [28]). Consequently, derivatives of order less than $\lfloor \frac{d}{2} \rfloor + 1$ of functions $h_n(f, \cdot)$ are uniformly bounded by some constant $\widetilde{C} > 0$. Therefore Property $(\mathcal{H}_0)$ is satisfied (by taking $L_{0,\alpha,\beta} = d^{d/2+2}\widetilde{C}$ and $a_{n,0,\alpha,\beta}(t) = 0$).

Let us now consider an integer $p \geq 0$. Let us suppose that $(\mathcal{H}_p)$ is satisfied and let us show that $(\mathcal{H}_{p+1})$ is then also satisfied. Let us notice that, since matrix $D(f)$ is nondegenerate, there exist two real numbers $c_0 > 0$ and $c_1 > 0$ such that, for every $u \in \mathbf{R}^d$, we have

$$c_0 |u|_\infty^2 \leq \langle u, D(f)u \rangle \leq c_1 |u|_\infty^2.$$

Let a real number $\beta > 0$ be fixed. There exists an integer $n_0 \geq 1$ such that, for all $u \in \mathbf{R}^d$ satisfying $|u|_\infty \leq n_0^{-\alpha}$, we have $\langle u, D(f)u \rangle < 1$ (e.g., any integer satisfying $n_0 > c_1^{1/(2\alpha)}$ is suitable). In the following, $n$ will be a nonnegative integer and $t$ a point in $\mathbf{R}^d$ satisfying $n \geq n_0$ and $|t|_\infty \leq n^{1/2-\alpha}$. We will then have $1 - \frac{\langle t, D(f)t \rangle}{2n} > \frac{1}{2} > 0$. The notation $O$ will only depend on $p$, $\alpha$, $\beta$ and $f$; for example, the notation $g_{n,t} = O(k_{n,t})$ means that there exists a real number $C > 0$ such that, for any integer $n \geq 1$ and any $t \in \mathbf{R}^d$ satisfying $|t|_\infty \leq n^{1/2-\alpha}$, we have $|g_{n,t}| \leq C \cdot |k_{n,t}|$. We will split $h_n(f,t)$ in pieces that we will estimate separately: $h_n(f,t) = \sum_{i=1}^{5} H_i(t,n)$.

*Part* 1. We start by estimating the following quantity:

$$(5) \qquad H_1(t,n) := \left(1 - \frac{\langle t, D(f)t \rangle}{2n}\right)^n - \exp\left\{-\frac{\langle t, D(f)t \rangle}{2}\right\}.$$

We will show that we have

$$(6) \quad \sum_{k=0}^{\lfloor d/2 \rfloor + 1} \sum_{1 \leq j_1, \ldots, j_k \leq d} \left|\frac{\partial^k}{\partial t_{j_1} \cdots \partial t_{j_k}} H_1(t,n)\right|$$
$$= O\left(\frac{1 + |t|_\infty^{(d+10)/2}}{\sqrt{n}} \exp\left\{-\frac{1}{2}\langle t, D(f)t \rangle \left(1 - \frac{d+4}{2n}\right)\right\}\right).$$

This term will contribute to the $a_{n,p+1,\alpha,\beta}$ term in (4) (for $p+1$ instead of $p$). Let us notice that we have

$$|H_1(t,n)| \leq c_1^2 \frac{|t|_\infty^4}{8n} \exp\left\{-\frac{1}{2}\langle t, D(f)t \rangle \left(1 - \frac{1}{n}\right)\right\}.$$

Let us now fix an integer $k \in \{1, \ldots, \lfloor \frac{d}{2} \rfloor + 1\}$ and $k$ indices $j_1, \ldots, j_k$ belonging to $\{1, \ldots, d\}$. In the following, we will denote by $Q_k$ the set of partitions $\mathcal{A} = \{\mathcal{A}_1, \ldots, \mathcal{A}_m\}$ of $\{1, \ldots, k\}$ in nonempty subsets. Let us notice that, for any $C^k$-regular function $b : \mathbf{R}^d \to \mathbf{R}$, we have

$$\frac{\partial^k}{\partial t_{j_1} \cdots \partial t_{j_k}}\left(\left(b\left(\frac{t}{\sqrt{n}}\right)\right)^n\right) = \sum_{\mathcal{A} = \{\mathcal{A}_1, \ldots, \mathcal{A}_m\} \in Q_k} g_n(\mathcal{A}, b)(t),$$

with

$$(7) \quad g_n(\mathcal{A}, b)(t) := n(n-1)\cdots(n-m+1)\left(b\left(\frac{t}{\sqrt{n}}\right)\right)^{n-m}$$
$$\times \prod_{p=1}^{m}\left(\frac{\partial^{\#\mathcal{A}_p} b}{\partial t_{j_{l_1^{(p)}}} \cdots \partial t_{j_{l_{\#\mathcal{A}_p}^{(p)}}}}\right)\left(\frac{t}{\sqrt{n}}\right) \frac{1}{n^{k/2}},$$

if $\mathcal{A} = \{\mathcal{A}_1, \ldots, \mathcal{A}_m\}$ with $\mathcal{A}_p := \{l_1^{(p)}, \ldots, l_{\#\mathcal{A}_p}^{(p)}\}$. In the following, we will consider that $b = 1 - \frac{1}{2}\langle \cdot, D(f) \cdot \rangle$ or $b = \exp\{-\frac{1}{2}\langle \cdot, D(f) \cdot \rangle\}$. Let $\mathcal{A} = \{\mathcal{A}_1, \ldots, \mathcal{A}_m\} \in$ $Q_k$. We denote by $m_0(\mathcal{A})$ the number of $\mathcal{A}_i \in \mathcal{A}$ which contains only one point. Then, we have $2m \leq m_0(\mathcal{A}) + k$. Indeed, we have

$$k = \sum_{p=1}^{m} \#\mathcal{A}_p \geq m_0(\mathcal{A}) + 2(m - m_0(\mathcal{A})) = 2m - m_0(\mathcal{A}).$$

(i) Let us suppose that $2m < m_0(\mathcal{A}) + k$. Using the fact that $(\frac{\partial}{\partial t_j} b)(\frac{t}{\sqrt{n}}) = O(\frac{|t|_\infty}{\sqrt{n}})$ and that the derivatives of order at least 2 of $b$ taken in $\frac{t}{\sqrt{n}}$ are bounded, we establish the following estimate:

$$|g_n(\mathcal{A}, b)(t)| \leq n^m \exp\left\{-\frac{n-m}{2n}\langle t, D(f)t\rangle\right\} O\left(\left(\frac{|t|_\infty}{\sqrt{n}}\right)^{m_0(\mathcal{A})}\right) n^{-k/2}$$

$$= O\left(n^{m-(m_0(\mathcal{A})+k)/2} |t|_\infty^{m_0(\mathcal{A})} \exp\left\{-\frac{1}{2}\langle t, D(f)t\rangle \left(1 - \frac{m}{n}\right)\right\}\right)$$

$$= O\left(n^{-1/2} |t|_\infty^{m_0(\mathcal{A})} \exp\left\{-\frac{1}{2}\langle t, D(f)t\rangle \left(1 - \frac{m}{n}\right)\right\}\right).$$

(ii) Let us suppose now that $2m = m_0(\mathcal{A}) + k$. Then each $\mathcal{A}_i$ contains at most two points and we show that we have

$$\left|g_n\left(\mathcal{A}, 1 - \frac{\langle \cdot, D(f) \cdot \rangle}{2}\right)(t) - g_n\left(\mathcal{A}, \exp\left\{-\frac{1}{2}\langle \cdot, D(f) \cdot \rangle\right\}\right)(t)\right|$$

$$= O\left(\frac{1 + |t|_\infty^{m_0(\mathcal{A})+4}}{n} \exp\left\{-\frac{1}{2}\langle t, D(f)t\rangle \left(1 - \frac{m+1}{n}\right)\right\}\right).$$

Effectively, let us notice that we have

$$\frac{\partial}{\partial t_j}\left(1 - \frac{\langle \cdot, D(f) \cdot \rangle}{2} - \exp\left\{-\frac{1}{2}\langle \cdot, D(f) \cdot \rangle\right\}\right)\left(\frac{t}{\sqrt{n}}\right)$$

$$= \sum_{l=1}^{d} D(f)_{j,l} \frac{t_l}{\sqrt{n}} \left(\exp\left\{-\frac{1}{2n}\langle t, D(f)t\rangle\right\} - 1\right) = O\left(\frac{|t|_\infty^3}{n\sqrt{n}}\right),$$

and we have

$$\left(1 - \frac{\langle t, D(f)t\rangle}{2n}\right)^{n-m} - \exp\left\{-\frac{n-m}{2n}\langle t, D(f)t\rangle\right\}$$

$$= O\left(\frac{|t|_\infty^4}{n} \exp\left\{-\frac{1}{2}\langle t, D(f)t\rangle\left(1 - \frac{m+1}{n}\right)\right\}\right),$$

by using formulae $|a^{n-m} - b^{n-m}| \leq (n-m)\max(|a|, |b|)^{n-m-1}|a-b|$ and $|e^{-u} - 1 - u| \leq \frac{u^2}{2}$. Moreover, we have

$$\frac{\partial^2}{\partial t_j \partial t_{j'}}\left(1 - \frac{\langle \cdot, D(f) \cdot \rangle}{2} - \exp\left\{-\frac{1}{2}\langle \cdot, D(f) \cdot \rangle\right\}\right)\left(\frac{t}{\sqrt{n}}\right)$$

$$= D(f)_{j,j'}\left(\exp\left\{-\frac{1}{2n}\langle t, D(f)t\rangle\right\} - 1\right)$$

$$- \sum_{l,m=1}^{d} D(f)_{j,l}\frac{t_l}{\sqrt{n}}D(f)_{j',m}\frac{t_m}{\sqrt{n}}\exp\left\{-\frac{1}{2n}\langle t, D(f)t\rangle\right\} = O\left(\frac{|t|_\infty^2}{n}\right).$$

Moreover, we recall that, for $b = 1 - \frac{1}{2}\langle\cdot, D(f)\cdot\rangle$ or $b = \exp\{-(1/2)\langle\cdot, D(f)\cdot\rangle\}$, we have $(\frac{\partial}{\partial t_j}b)(\frac{t}{\sqrt{n}}) = O(\frac{|t|_\infty}{\sqrt{n}})$ and that the derivatives of order at least 2 of $b$ taken in $\frac{t}{\sqrt{n}}$ are bounded.

Therefore, according to (7), the previous estimates and

$$\prod_{i=0}^{m} a_i - \prod_{i=0}^{m} b_i = \sum_{j=0}^{m}\left(\prod_{k=0}^{j-1} b_k\right)(a_j - b_j)\left(\prod_{l=j+1}^{m} a_l\right),$$

we get

$$\left|g_n\left(\mathcal{A}, 1 - \frac{\langle\cdot, D(f)\cdot\rangle}{2}\right)(t) - g_n\left(\mathcal{A}, \exp\left\{-\frac{1}{2}\langle\cdot, D(f)\cdot\rangle\right\}\right)(t)\right|$$

$$= O\left(n^m\left[\left(\frac{|t|_\infty}{\sqrt{n}}\right)^{m_0(\mathcal{A})+2} + \frac{|t|_\infty^{m_0(\mathcal{A})+4}}{\sqrt{n}^{m_0(\mathcal{A})+2}}\right]\right.$$

$$\left. \times \exp\left\{-\frac{1}{2}\langle t, D(f)t\rangle\left(1 - \frac{m+1}{n}\right)\right\}\frac{1}{n^{k/2}}\right)$$

$$= O\left(n^m\frac{1 + |t|_\infty^{m_0(\mathcal{A})+4}}{\sqrt{n}^{m_0(\mathcal{A})+2}}\exp\left\{-\frac{1}{2}\langle t, D(f)t\rangle\left(1 - \frac{m+1}{n}\right)\right\}\frac{1}{n^{k/2}}\right)$$

and we have $2m = m_0(\mathcal{A}) + k$.

*Part 2.* Hence, we have to study the quantity

$$D_n(t) := \mathbf{E}_\nu\left[\exp\left\{\frac{i\langle t, S_n(f)\rangle}{\sqrt{n}}\right\}\right] - \left(1 - \frac{\langle t, D(f)t\rangle}{2n}\right)^n,$$

which we split as follows:

$$(8)\quad D_n(t) = \sum_{l=0}^{n-1}\left(1 - \frac{\langle t, D(f)t\rangle}{2n}\right)^l \mathbf{E}_\nu\left[Y \circ T^l \cdot \exp\left\{\frac{i\langle t, S_{n-(l+1)}(f)\rangle}{\sqrt{n}}\right\} \circ T^{l+1}\right],$$

with

$$Y := \exp\left\{\frac{i\langle t, f\rangle}{\sqrt{n}}\right\} - \left(1 - \frac{\langle t, D(f)t\rangle}{2n}\right).$$

*Part 3.* Let us fix $M := p + 3$. Let us consider the nonnegative integers $a_1(n), \ldots, a_M(n)$ given by the formulae

$$a_1 := \left\lceil -\frac{\ln(n)}{\ln(\delta_{r_0})}\right\rceil,$$

$$a_j := \left\lceil (r_0-1)(a_1+\cdots+a_{j-1}) - \frac{\ln(n^{(d+5+\beta)/2}P_{r_0}(n))}{\ln(\delta_{r_0})} \right\rceil,$$

where $P_{r_0}$ and $\delta_{r_0}$ are, respectively, a polynomial function and a real number as in Property $(\mathcal{P}_{r_0})$. Let us write $A_0 := 0$ and $A_j := \sum_{k=1}^{j} a_j$. We notice that there exists a real number $\kappa > 0$ such that, for any integer $n \geq 1$ and any $j = 1, \ldots, M$, we have $a_j < \frac{\kappa n^{\alpha/2}}{M}$. Therefore, we have $a_1 + \cdots + a_M = O(n^\theta)$ for any real number $\theta > 0$.

*Part 4.* Let us define

(9)
$$H_2(t,n) := \sum_{n-\kappa n^{\alpha/2} \leq l \leq n-1} \left(1 - \frac{\langle t, D(f)t\rangle}{2n}\right)^l$$
$$\times \mathbf{E}_\nu\left[Y \cdot \exp\left\{\frac{i\langle t, S_{n-(l+1)}(f)\rangle}{\sqrt{n}}\right\} \circ T\right].$$

Let us prove that we have

(10)
$$\sum_{k=0}^{\lfloor d/2 \rfloor + 1} \sum_{j_1,\ldots,j_k=1,\ldots,d} \left|\frac{\partial^k}{\partial t_{j_1}\cdots\partial t_{j_k}} H_2(t,n)\right|$$
$$= O\left(n^{\alpha/2}\frac{(1+|t|_\infty^{(d+4)/2})}{\sqrt{n}}\right.$$
$$\left.\times \exp\left\{-\frac{\langle t, D(f)t\rangle}{2}\left(1 - \frac{\kappa}{n^{1-\alpha/2}} - \frac{d+2}{2n}\right)\right\}\right).$$

This term will contribute to the $a_{n,p+1,\alpha,\beta}$ term in (4) (for $p+1$ instead of $p$). Let us consider an integer $l$ satisfying $n - \kappa n^{\alpha/2} \leq l \leq n-1$ and an integer $k \geq 0$ and $k$ indices $j_1, \ldots, j_k$ in $\{1, \ldots, d\}$. First, let us notice that we have

$$\frac{\partial^k}{\partial t_{j_1}\cdots\partial t_{j_k}} \exp\left\{\frac{i\langle t, S_{n-(l+1)}(f)\rangle}{\sqrt{n}}\right\}$$
$$= \frac{i^k \prod_{p=1}^k S_{n-(l+1)}(f_{j_p})}{n^{k/2}} \exp\left\{\frac{i\langle t, S_{n-(l+1)}(f)\rangle}{\sqrt{n}}\right\} = O(1),$$

since we have $n - (l+1) \leq \kappa n^{\alpha/2}$ and $0 < \alpha < \frac{1}{2}$. Second, we have

$$\frac{\partial^k}{\partial t_{j_1}\cdots\partial t_{j_k}} Y = O\left(\frac{1+|t|_\infty}{\sqrt{n}}\right).$$

Indeed, $Y$ is in $O(\frac{|t|_\infty}{\sqrt{n}})$ and derivatives of $t \mapsto \exp\{\frac{i\langle t,f\rangle}{\sqrt{n}}\}$ are in $O(\frac{1}{\sqrt{n}})$. Moreover, derivatives of first order of $t \mapsto 1 - \frac{\langle t,D(f)t\rangle}{2n}$ are in $O(\frac{|t|_\infty}{n})$, its derivatives of order 2 are in $O(\frac{1}{n})$, its derivatives of order at least 3 are null. Now, let us show that we have

$$\frac{\partial^k}{\partial t_{j_1} \cdots \partial t_{j_k}} \left(1 - \frac{\langle t, D(f)t \rangle}{2n}\right)^l = O(b_{n,l}(t)), \tag{11}$$

with

$$b_{n,l}(t) := \sum_{m=0}^{\min(\lfloor d/2 \rfloor + 1, l)} \left(l \cdots (l-m+1) \frac{1+|t|_\infty^m}{n^m}\right)$$

$$\times \exp\left\{-\frac{1}{2n}\langle t, D(f)t \rangle (l-m)\right\}, \tag{12}$$

[with the convention $l \cdots (l - m + 1) = 1$ if $m = 0$]. Estimation (11) holds for $k = 0$ (since $|1 - u| \leq e^{-u}$ for any real number $u \in [0; 1]$). Let us suppose now $k \geq 1$. Since derivatives of order at least 3 of $t \mapsto 1 - \frac{\langle t, D(f)t \rangle}{2n}$ are null, we have

$$\frac{\partial^k}{\partial t_{j_1} \cdots \partial t_{j_k}} \left(1 - \frac{\langle t, D(f)t \rangle}{2n}\right)^l$$

$$= \sum_{\mathcal{A} \in \mathcal{B}_k} l \cdots (l - m + 1) \left(1 - \frac{\langle t, D(f)t \rangle}{2n}\right)^{l-m}$$

$$\times \prod_{p=1}^m \left(\frac{\partial^{\#\mathcal{A}_p}(1 - \langle \cdot, D(f) \cdot \rangle/2)}{\partial t_{j_1^{(p)}} \cdots \partial t_{j_{\#\mathcal{A}_p}^{(p)}}}\right) \left(\frac{t}{\sqrt{n}}\right) \frac{1}{n^{k/2}},$$

where we denote by $\mathcal{B}_k$ the set of partitions $\mathcal{A} = \{\mathcal{A}_1, \ldots, \mathcal{A}_m\}$ of $\{1, \ldots, k\}$ in subsets of at most two points. Let us consider such a partition $\mathcal{A} = \{\mathcal{A}_1, \ldots, \mathcal{A}_m\} \in \mathcal{B}_k$. If $m \geq l + 1$, then we have

$$l \cdots (l - m + 1)\left(1 - \frac{\langle t, D(f)t \rangle}{2n}\right)^{l-m}$$

$$\times \prod_{p=1}^m \left(\frac{\partial^{\#\mathcal{A}_p}(1 - \langle \cdot, D(f) \cdot \rangle/2)}{\partial t_{j_1^{(p)}} \cdots \partial t_{j_{\#\mathcal{A}_p}^{(p)}}}\right) \left(\frac{t}{\sqrt{n}}\right) \frac{1}{n^{k/2}} = 0.$$

Let us suppose now that $m \leq l$. Since we have $2m = m_0(\mathcal{A}) + k$, we get

$$\left| \frac{l!}{(l-m)!} \left(1 - \frac{\langle t, D(f)t \rangle}{2n}\right)^{l-m} \prod_{p=1}^m \left(\frac{\partial^{\#\mathcal{A}_p}(1 - \langle \cdot, D(f) \cdot \rangle/2)}{\partial t_{j_1^{(p)}} \cdots \partial t_{j_{\#\mathcal{A}_p}^{(p)}}}\right) \left(\frac{t}{\sqrt{n}}\right) \frac{1}{n^{k/2}} \right|$$

$$\leq \frac{l!}{(l-m)!} \left(1 - \frac{\langle t, D(f)t \rangle}{2n}\right)^{l-m}$$

$$\times \left(\frac{d\sup_{j,j'}|D(f)_{j,j'}|\cdot |t|_\infty}{\sqrt{n}}\right)^{m_0(\mathcal{A})} \left(\sup_{j,j'}|D(f)_{j,j'}|\right)^{m-m_0(\mathcal{A})} \frac{1}{n^{k/2}}$$

$$\leq \frac{l!}{(l-m)!}\left(1-\frac{\langle t,D(f)t\rangle}{2n}\right)^{l-m}\left(\frac{|t|_\infty^{m_0(\mathcal{A})}}{n^m}\right)\left(1+d\sup_{j,j'}|D(f)_{j,j'}|\right)^{d/2+1}.$$

Hence, we have proved (11). Let us prove now that we have

(13) $$\sum_{l=0}^{n-1} b_{n,l}(t) = O\left(\min\left(n, \frac{n}{|t|_\infty^2}\right)\right).$$

We have
$$b_{n,l}(t) = \sum_{m=0}^{\min(\lfloor d/2\rfloor+1, l)} b_{n,l,m}(t),$$

with
$$b_{n,l,m}(t) := l\cdots(l-m+1)\frac{1+|t|_\infty^m}{n^m}\exp\left\{-\frac{1}{2n}\langle t, D(f)t\rangle(l-m)\right\}.$$

We have
$$\sum_{l=0}^{n-1} b_{n,l}(t) = \sum_{m=0}^{\lfloor d/2\rfloor+1}\sum_{l=m}^{n-1} b_{n,l,m}(t).$$

Let us consider an integer $m \leq \lfloor \frac{d}{2}\rfloor + 1$. If $|t|_\infty \leq 1$, then we have

$$\sum_{l=m}^{n-1} b_{n,l,m}(t) \leq 2\sum_{l=m}^{n-1}\frac{l^m}{n^m}\exp\left\{-\frac{1}{2n}\langle t,D(f)t\rangle(l-m)\right\} \leq 2n.$$

If $|t|_\infty > 1$, then we have

$$\sum_{l=m}^{n-1} b_{n,l,m}(t) \leq 2\sum_{l=m}^{n-1} l\cdots(l-m+1)\frac{|t|_\infty^m}{n^m}\exp\left\{-\frac{1}{2n}\langle t,D(f)t\rangle(l-m)\right\}$$

$$\leq 2\frac{|t|_\infty^m}{n^m}\sum_{l\geq m} l\cdots(l-m+1)\exp\left\{-\frac{1}{2n}c_0|t|_\infty^2(l-m)\right\}$$

$$\leq 2\frac{|t|_\infty^m}{n^m}\frac{m!}{(1-\exp\{-(1/2n)c_0|t|_\infty^2\})^{m+1}}$$

$$\leq 2\frac{|t|_\infty^m}{n^m}\frac{m!}{(\exp\{-c_0/2\}(1/2n)c_0|t|_\infty^2)^{m+1}}$$

$$\leq 2\frac{|t|_\infty^m}{n^m}\frac{m!(2n)^{m+1}}{(\exp\{-c_0/2\}c_0|t|_\infty^2)^{m+1}}$$

$$\leq O\left(\frac{n}{|t|_\infty^{m+2}}\right) = O\left(\frac{n}{|t|_\infty^2}\right).$$

*Part* 5. For each nonnegative integer $l$ satisfying $n - (l+1) \geq \lceil \kappa n^{\alpha/2} \rceil$, we use the following decomposition of $S_{n-(l+1)}(f)$:

$$(14) \quad S_{n-(l+1)}(f) = \left( \sum_{j=1}^{M} S_{a_j}(f) \circ T^{A_{j-1}} \right) + S_{M_{n,l}}(f) \circ T^{A_M},$$

with $M_{n,l} := n - (l+1) - A_M$. Let us define

$$F_j^{(l)} := \exp\left\{ \frac{i\langle t, S_{a_j}(f) \rangle}{\sqrt{n}} \right\} \quad \text{and} \quad G^{(l)} := \exp\left\{ \frac{i\langle t, S_{M_{n,l}}(f) \rangle}{\sqrt{n}} \right\}.$$

We have

$$\mathbf{E}_\nu \left[ Y \circ T^l \cdot \exp\left\{ \frac{i\langle t, S_{n-(l+1)}(f) \rangle}{\sqrt{n}} \right\} \circ T^{l+1} \right]$$

$$= \mathbf{E}_\nu \left[ Y \left( \prod_{j=1}^{M} F_j^{(l)} \circ T^{1+A_{j-1}} \right) G^{(l)} \circ T^{1+A_M} \right].$$

We start by estimating the following quantity:

$$h_3(t, n, l) := \mathbf{E}_\nu \left[ Y \cdot F_1^{(l)} \circ T \left( \prod_{j=2}^{M} (F_j^{(l)} \circ T^{1+A_{j-1}} - 1) \right) G^{(l)} \circ T^{1+A_M} \right].$$

Let us show that we have

$$\sum_{k=0}^{\lfloor d/2 \rfloor + 1} \sum_{j_1,\ldots,j_k=1,\ldots,d} \left| \frac{\partial^k}{\partial t_{j_1} \cdots \partial t_{j_k}} h_3(t,n,l) \right| = O\left( \frac{1+|t|_\infty}{\sqrt{n}} \left( \frac{1+|t|_\infty}{n^{(1-\alpha)/2}} \right)^{M-1} \right)$$

$$= O\left( \frac{1+|t|_\infty^{p+3}}{\sqrt{n} \cdot n^{((1-\alpha)/2)(p+2)}} \right).$$

Effectively, for all $k = 0, \ldots, \lfloor \frac{d}{2} \rfloor + 1$ and all indices $j_1, \ldots, j_k \in \{1, \ldots, d\}$, we have

$$\frac{\partial^k}{\partial t_{j_1} \cdots \partial t_{j_k}} Y = O\left( \frac{1+|t|_\infty}{\sqrt{n}} \right) \quad \text{and} \quad \frac{\partial^k}{\partial t_{j_1} \cdots \partial t_{j_k}} F_1^{(l)} = O(1)$$

and

$$F_j^{(l)} - 1 = O\left( \frac{|t|_\infty}{n^{(1-\alpha)/2}} \right) \quad \text{and} \quad \frac{\partial^k}{\partial t_{j_1} \cdots \partial t_{j_k}} (F_j^{(l)} - 1) = O\left( \frac{1}{n^{(1-\alpha)/2}} \right)$$

and

$$\left\| \frac{\partial^k}{\partial t_{j_1} \cdots \partial t_{j_k}} G^{(l)} \right\|_{L^1} = O(1),$$

by using $|e^{iu} - 1| \leq |u|$ and $a_j < \kappa n^{\alpha/2}$ and the fact that $(\frac{S_n(f)}{\sqrt{n}})_{n \geq 1}$ is uniformly bounded in $L^p$, for all $p \in [1, +\infty[$. Let us define

$$H_3(t,n) := \sum_{l=0}^{n-\lfloor \kappa n^{\alpha/2} \rfloor - 1} \left(1 - \frac{\langle t, D(f)t \rangle}{2n}\right)^l h_3(t,n,l).$$

According to (11) and (13) we have

(15)
$$\sum_{k=0}^{\lfloor d/2 \rfloor + 1} \sum_{j_1, \ldots, j_k = 1, \ldots, d} \left| \frac{\partial^k}{\partial t_{j_1} \cdots \partial t_{j_k}} H_3(t,n) \right| = O\left(\frac{1 + |t|_\infty^{p+1}}{n^{(p+1-\alpha(p+2))/2}}\right)$$
$$= O\left(\frac{1 + |t|_\infty^{p+1}}{n^{(1/2-\alpha)(p+1)}}\right).$$

This term will contribute to the first term of estimate (4) (for $p + 1$ instead of $p$).

*Part* 6. (Heart of the proof.) It remains to estimate the following term:

$$\sum_{l=0}^{n-\lfloor \kappa n^{\alpha/2} \rfloor - 1} \left(1 - \frac{\langle t, D(f)t \rangle}{2n}\right)^l$$
$$\times \sum_{\varepsilon = (\varepsilon_1, \ldots, \varepsilon_m)} \mathbf{E}_\nu \left[ Y \left( \prod_{j=1}^M \varepsilon_j \circ T^{1+A_{j-1}} \right) G^{(l)} \circ T^{1+A_M} \right],$$

the second sum being taken over the set of $\varepsilon = (\varepsilon_1, \ldots, \varepsilon_M) \in \prod_{j=1}^M \{-1; F_j^{(l)}\}$, with $\varepsilon_1 := F_1^{(l)}$ and with at least one $\varepsilon_j$ equal to $-1$. Let an integer $l = 0, \ldots, n - \lceil \kappa n^{\alpha/2} \rceil - 1$ and such a vector $\varepsilon = (\varepsilon_1, \ldots, \varepsilon_M)$ be given. We define $j_0 := \max\{j \geq 2 : \varepsilon_j = -1\}$. Then we define

$$D_{l,\varepsilon}(n,t) := Y \prod_{j=1}^{j_0 - 1} \varepsilon_j \circ T^{1+A_{j-1}}$$

and

$$E_{l,\varepsilon}(n,t) := \exp\left\{ \frac{i}{\sqrt{n}} \langle t, S_{n-(l+1)-A_{j_0}}(f) \rangle \right\}.$$

Therefore, we have

$$\mathbf{E}_\nu \left[ Y \left( \prod_{j=1}^M \varepsilon_j \circ T^{1+A_{j-1}} \right) G^{(l)} \circ T^{1+A_M} \right] = -\mathbf{E}_\nu [D_{l,\varepsilon}(n,t) \cdot E_{l,\varepsilon}(n,t) \circ T^{1+A_{j_0}}].$$

*First step*: *control of* $\operatorname{Cov}_\nu(D_{l,\varepsilon}(n,t), E_{l,\varepsilon}(n,t) \circ T^{1+A_{j_0}})$. Let us prove that, for all $k = 0, \ldots, \lfloor \frac{d}{2} \rfloor + 1$ and all $j_1, \ldots, j_k \in \{1, \ldots, d\}$, we have

$$(16) \quad \left| \frac{\partial^k}{\partial t_{j_1} \cdots \partial t_{j_k}} \operatorname{Cov}_\nu(D_{l,\varepsilon}(n,t), E_{l,\varepsilon}(n,t) \circ T^{1+A_{j_0}}) \right| = O\left( \frac{1 + |t|_\infty}{n^{(d+3+\beta)/2}} \right).$$

We will use Property $(\mathcal{P}_{r_0})$. Let us notice that the functions $D_{l,\varepsilon}(n,t)$ and $E_{l,\varepsilon}(n,t)$ are of the following form:

$$D_{l,\varepsilon}(n,t) = Y \prod_{j=1}^{A_{j_0}-1} \alpha_j \circ T^j \quad \text{and} \quad E_{l,\varepsilon}(n,t) = \prod_{j=0}^{n-(l+1)-A_{j_0}-1} \exp\left\{ \frac{i\langle t, f \rangle}{\sqrt{n}} \right\} \circ T^j$$

for some $\alpha_j \in \{1, -1, \exp\{\frac{i\langle t,f\rangle}{\sqrt{n}}\}\}$. First, let us explain how we get (16) when $k = 0$. Let us notice that $Y$ is in $O(\frac{|t|_\infty}{\sqrt{n}})$ and is $\eta$-Hölder continuous with Hölder constant in $O(\frac{|t|_\infty}{\sqrt{n}})$. Moreover, $\|\alpha_j\|_\infty = 1$ and $\alpha_j$ are $\eta$-Hölder continuous with Hölder constant uniformly bounded in $O(\frac{|t|_\infty}{\sqrt{n}})$. Therefore, according to Property $(\mathcal{P}_{r_0})$, we get

$$|\operatorname{Cov}_\nu(D_{l,\varepsilon}(n,t), E_{l,\varepsilon}(n,t) \circ T^{1+A_{j_0}})|$$
$$\leq \left( \frac{|t|_\infty}{\sqrt{n}} + nO\left(\frac{|t|_\infty}{\sqrt{n}}\right) \right) P_{r_0}(n - (l+1) - A_{j_0}) \delta_{r_0}^{1+A_{j_0} - r_0 A_{j_0-1}}$$
$$\leq O(|t|_\infty \sqrt{n}) P_{r_0}(n) \delta_{r_0}^{1+A_{j_0} - r_0 A_{j_0-1}}$$
$$\leq O(|t|_\infty \sqrt{n}) \frac{1}{n^{(d+4+\beta)/2}},$$

according to the fact that

$$P_{r_0}(n) \delta_{r_0}^{1+A_{j_0} - r_0 A_{j_0-1}} \leq \frac{1}{n^{(d+4+\beta)/2}},$$

(see the definition of $a_{j_0}$). Let us suppose now $k \geq 1$. The partial derivatives of $Y$ relative to $t$ are in $O(\frac{1+|t|_\infty}{\sqrt{n}})$ and are $\eta$-Hölder continuous with Hölder constant in $O(\frac{1+|t|_\infty}{\sqrt{n}})$. Moreover, the partial derivatives of $\alpha_j$ relative to $t$ are uniformly bounded in $O(\frac{1}{\sqrt{n}})$ and are $\eta$-Hölder continuous with Hölder constant in $O(\frac{1+|t|_\infty}{\sqrt{n}})$. Therefore, the derivative of order $k' \geq 1$ of $\prod_{j=1}^{A_{j_0}-1} \alpha_j \circ T^j$ is a sum of $(A_{j_0-1})^{k'}$ terms of the following form: $\prod_{j=1}^{A_{j_0}-1} \beta_j \circ T^j$, where $\beta_j$ is equal to $\alpha_j$ or to some derivative of $\alpha_j$ and with at least one $\beta_j$ equal to some derivative of $\alpha_j$. Therefore, according to Property $(\mathcal{P}_{r_0})$, for all integers $k_1 \geq 1$ and $k_2 \geq 1$ such that $k_1 + k_2 = k$, and all $i_1, \ldots, i_{k_1}, j_1, \ldots, j_{k_2}$ in $\{1,\ldots,d\}$, we have

$$\left|\mathrm{Cov}_\nu\left(\frac{\partial^{k_1}}{\partial t_{i_1}\cdots\partial t_{i_{k_1}}}D_{l,\varepsilon}(n,t),\frac{\partial^{k_2}}{\partial t_{j_1}\cdots\partial t_{j_{k_2}}}E_{l,\varepsilon}(n,t)\circ T^{1+A_{j_0}}\right)\right|$$

$$\leq A_{j_0-1}{}^k O\left(\frac{1+|t|_\infty}{\sqrt{n}}+n\frac{1+|t|_\infty}{\sqrt{n}}\right)P_{r_0}(n)\delta_{r_0}{}^{1+A_{j_0}-r_0 A_{j_0-1}}$$

$$\leq O((1+|t|_\infty)\sqrt{n})A_{j_0-1}{}^{d/2+1}P_{r_0}(n)\delta_{r_0}{}^{1+A_{j_0}-r_0 A_{j_0-1}}.$$

We conclude by using the facts that $P_{r_0}(n)\delta_{r_0}{}^{1+A_{j_0}-r_0 A_{j_0-1}}\leq\frac{1}{n^{(d+5+\beta)/2}}$ and that $A_{j_0-1}$ is in $O(\ln(n))$ (see the definition of $a_j$).

We define

$$H_4(t,n):=\sum_{l=0}^{n-\lfloor\kappa n^{\alpha/2}\rfloor-1}\left(1-\frac{\langle t,D(f)t\rangle}{2n}\right)^l$$

(17)

$$\times\sum_{\varepsilon=(\varepsilon_1,\ldots,\varepsilon_M)}\mathrm{Cov}_\nu(D_{l,\varepsilon}(n,t),E_{l,\varepsilon}(n,t)\circ T^{1+A_{j_0}}).$$

According to the preceding and (11) and (13) we have

(18) $$\sum_{k=0}^{\lfloor d/2\rfloor+1}\sum_{j_1,\ldots,j_k\in\{1,\ldots,d\}}\left|\frac{\partial^k}{\partial t_{j_1}\cdots\partial t_{j_k}}H_4(t,n)\right|=O\left(\frac{1}{n^{(d+1+\beta)/2}}\right).$$

This term will contribute to the $a_{n,p+1,\alpha,\beta}$ term in (4) (for $p+1$ instead of $p$). It remains to estimate the derivatives of the following quantity:

$$H_5(t,n):=\sum_{l=0}^{n-\lfloor\kappa n^{\alpha/2}\rfloor-1}\left(1-\frac{\langle t,D(f)t\rangle}{2n}\right)^l$$

(19)

$$\times\sum_{\varepsilon=(\varepsilon_1,\ldots,\varepsilon_M)}\mathbf{E}_\nu[D_{l,\varepsilon}(n,t)]\mathbf{E}_\nu[E_{l,\varepsilon}(n,t)\circ T^{1+A_{j_0}}].$$

*Second step*: control of the expectation of $D_{l,\varepsilon}(n,t)$. Let us show that we have

(20) $$\sup_{l=0,\ldots,n-\lfloor\kappa n^{\alpha/2}\rfloor-1}\sup_{\varepsilon=(\varepsilon_1,\ldots,\varepsilon_m)}\sum_{k=0}^{\lfloor d/2\rfloor+1}\sum_{j_1,\ldots,j_k\in\{1,\ldots,d\}}\left|\frac{\partial^k}{\partial t_{j_1}\cdots\partial t_{j_k}}\mathbf{E}_\nu[D_{l,\varepsilon}(n,t)]\right|$$
$$=O\left(\frac{1+|t|_\infty^3}{\sqrt{n}\cdot n^{1-\alpha}}\right).$$

Let us denote by $\mathcal{J}$ the following set:
$$\mathcal{J} := \{j = 1, \ldots, j_0 - 1 : \varepsilon_j = \mathcal{F}_j^{(l)}\}.$$
Let us recall that 1 belongs to $\mathcal{J}$. In the following, we denote $S_{\mathcal{J}}(g) := \sum_{j \in \mathcal{J}} S_{a_j}(g) \circ T^{1+A_{j-1}} = \sum_{j \in \mathcal{J}} \sum_{k=A_{j-1}+1}^{A_j} g \circ T^k$. We have

$$|\mathbf{E}_\nu[D_{l,\varepsilon}(n,t)]| = \left|\mathbf{E}_\nu\left[Y \cdot \exp\left\{\frac{i}{\sqrt{n}}\langle t, S_{\mathcal{J}}(f)\rangle\right\}\right]\right|$$

$$= \left|\mathbf{E}_\nu\left[\left(\exp\left\{\frac{i\langle t, f\rangle}{\sqrt{n}}\right\} - 1 + \frac{\langle t, D(f)t\rangle}{2n}\right)\exp\left\{\frac{i}{\sqrt{n}}\langle t, S_{\mathcal{J}}(f)\rangle\right\}\right]\right|.$$

*Case $k = 0$.* With the use of Taylor's formulae of order 2 and 1 for $e^{iu}$, we get

$$|\mathbf{E}_\nu[D_{l,\varepsilon}(n,t)]|$$
$$= \left|\mathbf{E}_\nu\left[\left(\frac{i\langle t,f\rangle}{\sqrt{n}} + \frac{1}{2n}(\langle t, D(f)t\rangle - \langle t,f\rangle^2)\right)\left(1 + \frac{i}{\sqrt{n}}\langle t, S_{\mathcal{J}}(f)\rangle\right)\right]\right|$$
$$+ O\left(\frac{|t|_\infty^3 n^\alpha}{n^{3/2}}\right)$$
$$= \left|\frac{1}{2n}\mathbf{E}_\nu[\langle t, D(f)t\rangle - \langle t,f\rangle^2 - 2\langle t,f\rangle\langle t, S_{\mathcal{J}}(f)\rangle]\right| + O\left(\frac{|t|_\infty^3 n^\alpha}{n^{3/2}}\right)$$
$$= \left|\frac{1}{2n}\langle t, (D(f) - \mathbf{E}_\nu[f^{\otimes 2}] - \mathbf{E}_\nu[f \otimes S_{\mathcal{J}}(f)] - \mathbf{E}_\nu[S_{\mathcal{J}}(f) \otimes f])t\rangle\right|$$
$$+ O\left(\frac{|t|_\infty^3 n^\alpha}{n^{3/2}}\right)$$
$$= O\left(\frac{|t|_\infty^2}{n^2}\right) + O\left(\frac{|t|_\infty^3 n^\alpha}{n^{3/2}}\right) = O\left(\frac{1+|t|_\infty^3}{\sqrt{n}n^{1-\alpha}}\right).$$

Term in $O(\frac{|t|_\infty^2}{n^2})$ comes from (2) and from the fact that $\mathbf{E}_\nu[f_j.f_{j'} \circ T^k]$ converges to 0 exponentially fast as $k$ goes to infinity [this is a consequence of Property $(\mathcal{P}_{r_0})$]. Effectively, since 1 is in $\mathcal{J}$, we have

$$(21) \qquad S_{\mathcal{J}}(f) = \sum_{k=1}^{a_1} f \circ T^k + \sum_{k' \geq a_1+1, k' \in \mathcal{L}} f \circ T^{k'},$$

for some set of integers $\mathcal{L}$, and we have $(\delta_{r_0})^{a_1} \leq \frac{1}{n}$.

*Case $k \geq 3$.* Let us recall that we have:
- $Y = \exp\{\frac{i\langle t,f\rangle}{\sqrt{n}}\} - 1 + \frac{1}{2n}\langle t, D(f)t\rangle = O(\frac{|t|_\infty}{\sqrt{n}})$;

- for any $j \in \{1, \ldots, d\}$, $\frac{\partial}{\partial t_j} Y = \frac{i f_j}{\sqrt{n}} \exp\{\frac{i \langle t, f \rangle}{\sqrt{n}}\} + \frac{1}{n}{}^T e_j \cdot D(f) \cdot t = O(\frac{1}{\sqrt{n}} + \frac{|t|_\infty}{n})$, where $e_j$ is the $j$th vector of the canonical basis of $\mathbf{R}^d$;
- for all $j, j' \in \{1, \ldots, d\}$, $\frac{\partial^2}{\partial t_j \partial t_{j'}} Y = -\frac{f_j f_{j'}}{n} \exp\{\frac{i \langle t, f \rangle}{\sqrt{n}}\} + \frac{1}{n} D(f)_{j,j'} = O(\frac{1}{n})$;
- for any integer $m \geq 3$ and any $(j_1, \ldots, j_m)$ in $\{1, \ldots, d\}^m$,

$$\frac{\partial^m}{\partial t_{j_1} \cdots \partial t_{j_m}} Y = i^m \frac{f_{j_1} \cdot \ldots \cdot f_{j_m}}{n^{m/2}} \exp\left\{\frac{i \langle t, f \rangle}{\sqrt{n}}\right\} = O\left(\frac{1}{n^{m/2}}\right);$$

- for any integer $m \geq 0$ and any $(j_1, \ldots, j_m)$ in $\{1, \ldots, d\}^m$,

$$\frac{\partial^m}{\partial t_{j_1} \cdots \partial t_{j_m}} \exp\left\{\frac{i}{\sqrt{n}} \langle t, S_{\mathcal{J}}(f) \rangle\right\}$$

$$= \left(\frac{i S_{\mathcal{J}}(f_{j_1})}{\sqrt{n}}\right) \cdot \ldots \cdot \left(\frac{i S_{\mathcal{J}}(f_{j_m})}{\sqrt{n}}\right) \exp\left\{\frac{i}{\sqrt{n}} \langle t, S_{\mathcal{J}}(f) \rangle\right\}$$

is in $O(\frac{1}{n^{m(1/2-\alpha/3)}})$, according to the fact that $a_j = O(n^{\alpha/3})$.

Hence, for any integer $k \geq 3$ and any $(j_1, \ldots, j_k)$ in $\{1, \ldots, d\}^k$, we have

$$\frac{\partial^k}{\partial t_{j_1} \cdots \partial t_{j_k}} \mathbf{E}_\nu[D_{l,\varepsilon}(n, t)] = O\left(\frac{1 + |t|_\infty}{\sqrt{n} n^{1-\alpha}}\right).$$

*Case $k = 1$.* Let $j \in \{1, \ldots, d\}$ be given. We have

$$\frac{\partial}{\partial t_j} \mathbf{E}_\nu[D_{l,\varepsilon}(n, t)]$$

$$= \mathbf{E}_\nu\left[\left(\frac{\partial}{\partial t_j} Y\right) \exp\left\{\frac{i}{\sqrt{n}} \langle t, S_{\mathcal{J}}(f) \rangle\right\}\right] + \mathbf{E}_\nu\left[Y \cdot \frac{\partial}{\partial t_j} \exp\left\{\frac{i}{\sqrt{n}} \langle t, S_{\mathcal{J}}(f) \rangle\right\}\right]$$

$$= \mathbf{E}_\nu\left[\left(\frac{i f_j}{\sqrt{n}} \exp\left\{\frac{i \langle t, f \rangle}{\sqrt{n}}\right\} + \frac{1}{n}{}^T e_j \cdot D(f) \cdot t\right) \exp\left\{\frac{i}{\sqrt{n}} \langle t, S_{\mathcal{J}}(f) \rangle\right\}\right]$$

$$+ \mathbf{E}_\nu\left[\left(\exp\left\{\frac{i \langle t, f \rangle}{\sqrt{n}}\right\} - 1 + \frac{1}{2n} \langle t, D(f) t \rangle\right) \frac{i S_{\mathcal{J}}(f_j)}{\sqrt{n}}\right.$$

$$\left.\times \exp\left\{\frac{i}{\sqrt{n}} \langle t, S_{\mathcal{J}}(f) \rangle\right\}\right]$$

$$= \mathbf{E}_\nu\left[\left(i \frac{f_j + S_{\mathcal{J}}(f_j)}{\sqrt{n}}\right) \exp\left\{\frac{i}{\sqrt{n}} \langle t, f + S_{\mathcal{J}}(f) \rangle\right\}\right]$$

$$+ \mathbf{E}_\nu\left[\left(-i \frac{S_{\mathcal{J}}(f_j)}{\sqrt{n}} + \frac{{}^T e_j \cdot D(f) \cdot t}{n}\right) \exp\left\{\frac{i}{\sqrt{n}} \langle t, S_{\mathcal{J}}(f) \rangle\right\}\right]$$

$$+ O\left(\frac{|t|_\infty^2}{\sqrt{n} n^{1-\alpha}}\right)$$

$$= \mathbf{E}_\nu\left[i\frac{f_j + S_\mathcal{J}(f_j)}{\sqrt{n}}\left(1 + \frac{i}{\sqrt{n}}\langle t, f + S_\mathcal{J}(f)\rangle\right)\right]$$

$$- \mathbf{E}_\nu\left[i\frac{S_\mathcal{J}(f_j)}{\sqrt{n}}\left(1 + \frac{i}{\sqrt{n}}\langle t, S_\mathcal{J}(f)\rangle\right)\right] + \frac{{}^T e_j \cdot D(f) \cdot t}{n} + O\left(\frac{|t|_\infty^2}{\sqrt{n}n^{1-\alpha}}\right)$$

$$= -\frac{1}{n}\mathbf{E}_\nu[(f_j + S_\mathcal{J}(f_j))\langle t, f + S_\mathcal{J}(f)\rangle]$$

$$+ \frac{1}{n}\mathbf{E}_\nu[S_\mathcal{J}(f_j)\langle t, S_\mathcal{J}(f)\rangle] + \frac{{}^T e_j \cdot D(f) \cdot t}{n} + O\left(\frac{|t|_\infty^2}{\sqrt{n}n^{1-\alpha}}\right)$$

$$= \frac{1}{n}{}^T e_j(D(f) - \mathbf{E}_\nu[f^{\otimes 2}] - \mathbf{E}_\nu[f \otimes S_\mathcal{J}(f)] - \mathbf{E}_\nu[S_\mathcal{J}(f) \otimes f])t$$

$$+ O\left(\frac{|t|_\infty^2}{\sqrt{n}n^{1-\alpha}}\right) = O\left(\frac{|t|_\infty}{n^2}\right) + O\left(\frac{|t|_\infty^2}{\sqrt{n}n^{1-\alpha}}\right) = O\left(\frac{1 + |t|_\infty^2}{\sqrt{n}n^{1-\alpha}}\right).$$

Term in $O(\frac{|t|_\infty}{n^2})$ comes from (2) and (21) and from the fact that $\mathbf{E}_\nu[f_j.f_{j'} \circ T^k]$ converges to 0 exponentially fast as $k$ goes to infinity.

*Case $k = 2$.* Let $j_1$ and $j_2$ be in $\{1, \ldots, d\}$. We have

$$\frac{\partial^2}{\partial t_{j_1} \partial t_{j_2}}\mathbf{E}_\nu[D_{l,\varepsilon}(n,t)]$$

$$= \mathbf{E}_\nu\left[Y\frac{\partial^2}{\partial t_{j_1} \partial t_{j_2}}\exp\left\{\frac{i}{\sqrt{n}}\langle t, S_\mathcal{J}(f)\rangle\right\}\right]$$

$$+ \mathbf{E}_\nu\left[\left(\frac{\partial}{\partial t_{j_1}}Y\right)\frac{\partial}{\partial t_{j_2}}\exp\left\{\frac{i}{\sqrt{n}}\langle t, S_\mathcal{J}(f)\rangle\right\}\right]$$

$$+ \mathbf{E}_\nu\left[\left(\frac{\partial}{\partial t_{j_2}}Y\right)\frac{\partial}{\partial t_{j_1}}\exp\left\{\frac{i}{\sqrt{n}}\langle t, S_\mathcal{J}(f)\rangle\right\}\right]$$

$$+ \mathbf{E}_\nu\left[\left(\frac{\partial^2}{\partial t_{j_1} \partial t_{j_2}}Y\right)\exp\left\{\frac{i}{\sqrt{n}}\langle t, S_\mathcal{J}(f)\rangle\right\}\right]$$

$$= \frac{1}{n}(-\mathbf{E}_\nu[f_{j_1}S_\mathcal{J}(f_{j_2})] - \mathbf{E}_\nu[f_{j_2}S_\mathcal{J}(f_{j_1})] - \mathbf{E}_\nu[f_{j_1}f_{j_2}] + D(f)_{j_1,j_2})$$

$$+ O\left(\frac{1 + |t|_\infty^2}{\sqrt{n}n^{1-\alpha}}\right)$$

$$= O\left(\frac{1}{n^2}\right) + O\left(\frac{1 + |t|_\infty^2}{\sqrt{n}n^{1-\alpha}}\right) = O\left(\frac{1 + |t|_\infty^2}{\sqrt{n}n^{1-\alpha}}\right).$$

*Third step*: *control of the expectation of* $E_{l,\varepsilon}(n,t)$. We define

$$n' = n'_{n,l,\varepsilon} := n - (l+1) - A_{j_0} \quad \text{and} \quad t' = t'_{n,l,\varepsilon} := t\sqrt{\frac{n'}{n}}.$$

We take $\beta' := \beta + d + 8$. According to the inductive hypothesis $(\mathcal{H}_p)$ applied to $(n', t')$, we have

$$\sum_{k=0}^{\lfloor d/2 \rfloor + 1} \sum_{j_1,\ldots,j_k=1,\ldots,d} \left| \frac{\partial^k}{\partial t_{j_1} \cdots \partial t_{j_k}} \left( \mathbf{E}_\nu[E_{l,\varepsilon}(n,t)] \right. \right.$$

$$\left. \left. - \exp\left\{ -\frac{\langle t, D(f)t \rangle}{2} \left(1 - \frac{l+1}{n} - \frac{A_{j_0}}{n}\right) \right\} \right) \right|$$

$$\leq L_{p,\alpha,\beta'} \frac{1 + |t'|_\infty^p}{n'^{p(1/2-\alpha)}} + a_{n',p,\alpha,\beta'}(t').$$

Hence, since $A_{j_0} \leq \kappa n^\alpha$, we have

$$\sum_{k=0}^{\lfloor d/2 \rfloor + 1} \sum_{j_1,\ldots,j_k=1,\ldots,d} \left| \frac{\partial^k}{\partial t_{j_1} \cdots \partial t_{j_k}} \mathbf{E}_\nu[E_{l,\varepsilon}(n,t)] \right|$$

(22)
$$\leq O\left( (1 + |t|_\infty^{d/2+1}) \exp\left\{ -\frac{\langle t, D(f)t \rangle}{2} \left(1 - \frac{l+1}{n} - \frac{\kappa}{n^{1-\alpha}}\right) \right\} \right)$$

$$+ L_{p,\alpha,\beta'} \frac{1 + |t'|_\infty^p}{n'^{p(1/2-\alpha)}} + a_{n',p,\alpha,\beta'}(t').$$

*Part* 7. (Conclusion.) To finish the proof of Proposition 2.7, we deduce from the preceding an estimate of the following quantity:

(23) $$\sum_{k=0}^{\lfloor d/2 \rfloor + 1} \sum_{j_1,\ldots,j_k=1,\ldots,d} \left| \frac{\partial^k}{\partial t_{j_1} \cdots \partial t_{j_k}} H_5(t,n) \right|,$$

where we denote by $H_5$ the quantity introduced in (19). According to (11), (20) and (22), we have

$$\sum_{k=0}^{\lfloor d/2 \rfloor + 1} \sum_{j_1,\ldots,j_k=1,\ldots,d} \left| \frac{\partial^k}{\partial t_{j_1} \cdots \partial t_{j_k}} H_5(t,n) \right|$$

$$= \sum_{l=0}^{n - \lfloor \kappa n^{\alpha/2} \rfloor - 1} O\left( b_{n,l}(t) \left( \frac{1 + |t|_\infty^3}{\sqrt{n} n^{1-\alpha}} \right) \right.$$

$$\left. \times \left( (1 + |t|_\infty^{d/2+1}) \exp\left\{ -\frac{\langle t, D(f)t \rangle}{2} \left(1 - \frac{l+1}{n} - \frac{\kappa}{n^{1-\alpha}}\right) \right\} \right. \right.$$

$$+ L_{p,\alpha,\beta'} \frac{1+|t'|_\infty^p}{n'^{p(1/2-\alpha)}} + a_{n',p,\alpha,\beta'}(t')\bigg)\bigg).$$

Let us now estimate each term of the right-hand side part of this inequality. We will use (13) in (b) and (c). In (a) and (d)–(f), we use the fact that $b_{n,l}(t)$ is in $O((1+|t|_\infty^{\lfloor d/2 \rfloor+1})\exp\{-\frac{1}{2n}\langle t, D(f)t\rangle(l-\frac{d}{2}-1)\})$.

(a) We have

$$\sum_{l=0}^{n-\lfloor \kappa n^{\alpha/2}\rfloor -1} b_{n,l}(t)\bigg(\frac{1+|t|_\infty^3}{\sqrt{n}n^{1-\alpha}}\bigg)(1+|t|_\infty^{d/2+1})$$

$$\times \exp\bigg\{-\frac{\langle t, D(f)t\rangle}{2}\bigg(1-\frac{l+1}{n}-\frac{\kappa}{n^{1-\alpha}}\bigg)\bigg\}$$

$$= O\bigg(\frac{1+|t|_\infty^{d+5}}{n^{1/2-\alpha}}\exp\bigg\{-\frac{\langle t, D(f)t\rangle}{2}\bigg(1-\frac{\kappa}{n^{1-\alpha}}-\frac{d+4}{2n}\bigg)\bigg\}\bigg).$$

(b) We have

$$\sum_{l=0}^{n-\lfloor \kappa n^{\alpha/2}\rfloor -1} b_{n,l}(t)\bigg(\frac{1+|t|_\infty^3}{\sqrt{n}n^{1-\alpha}}\bigg)\frac{|t'|_\infty^p}{n'^{p(1/2-\alpha)}}$$

$$\leq \sum_{l=0}^{n-\lfloor \kappa n^{\alpha/2}\rfloor -1} b_{n,l}(t)\bigg(\frac{1+|t|_\infty^3}{\sqrt{n}n^{1-\alpha}}\bigg)\frac{|t|_\infty^p}{n^{p(1/2-\alpha)}} \leq O\bigg(\frac{1+|t|_\infty^{p+1}}{n^{(p+1)(1/2-\alpha)}}\bigg).$$

(c) Let us notice that if $l \leq \lfloor \frac{n}{2}\rfloor - \lceil \kappa n^\alpha\rceil -1$, then we have $n' \geq \frac{n}{2}$, from which we get

$$\sum_{l=0}^{\lfloor n/2\rfloor-\lceil \kappa n^\alpha\rceil-1} b_{n,l}(t)\bigg(\frac{1+|t|_\infty^3}{\sqrt{n}n^{1-\alpha}}\bigg)\frac{1}{n'^{p(1/2-\alpha)}}$$

$$= O\bigg(\bigg(\frac{1+|t|_\infty^3}{\sqrt{n}n^{1-\alpha}}\bigg)\frac{1}{n^{p(1/2-\alpha)}}n\min\bigg(1,\frac{1}{|t|_\infty^2}\bigg)\bigg)$$

$$= O\bigg(\bigg(\frac{1+|t|_\infty}{n^{1/2-\alpha}}\bigg)\frac{1}{n^{p(1/2-\alpha)}}\bigg) = O\bigg(\frac{1+|t|_\infty^{p+1}}{n^{(p+1)(1/2-\alpha)}}\bigg).$$

(d) We have

$$\sum_{l=\lfloor n/2\rfloor-\lceil \kappa n^\alpha\rceil}^{n-\lfloor \kappa n^{\alpha/2}\rfloor-1} b_{n,l}(t)\bigg(\frac{1+|t|_\infty^3}{\sqrt{n}n^{1-\alpha}}\bigg)\frac{1}{n'^{p(1/2-\alpha)}}$$

$$= O\bigg(\frac{1+|t|_\infty^{\lfloor d/2\rfloor+4}}{\sqrt{n}n^{1-\alpha}}n\exp\bigg\{-\frac{\langle t, D(f)t\rangle}{2n}\bigg(\frac{n}{2}-\kappa n^\alpha - 2 - \frac{d}{2}-1\bigg)\bigg\}\bigg)$$

$$= O\left(\frac{1+|t|_\infty^{\lfloor d/2\rfloor+4}}{n^{1/2-\alpha}}\exp\left\{-\frac{\langle t, D(f)t\rangle}{2}\left(\frac{1}{2}-\frac{\kappa}{n^{1-\alpha}}-\frac{d+6}{2n}\right)\right\}\right).$$

(e) We have

$$\left(\int_{|t|_\infty\le n^{1/2-\alpha}}(1+|t|_\infty^\beta)\left(\sum_{l=0}^{\lfloor n/2\rfloor-\lceil\kappa n^\alpha\rceil-1}b_{n,l}(t)\left(\frac{1+|t|_\infty^3}{\sqrt{n}n^{1-\alpha}}\right)a_{n',p,\alpha,\beta'}(t')\right)^2 dt\right)^{1/2}$$

$$= O\left(\left(\int_{|t|_\infty\le n^{1/2-\alpha}}(1+|t|_\infty^\beta)\right.\right.$$
$$\left.\left.\times\left(\sum_{l=0}^{\lfloor n/2\rfloor-\lceil\kappa n^\alpha\rceil-1}\left(\frac{1+|t|_\infty^{\lfloor d/2\rfloor+4}}{\sqrt{n}n^{1-\alpha}}\right)a_{n',p,\alpha,\beta'}(t')\right)^2 dt\right)^{1/2}\right)$$

$$= O\left(\sum_{l=0}^{\lfloor n/2\rfloor-\lceil\kappa n^\alpha\rceil-1}\left(\int_{|t|_\infty\le n^{1/2-\alpha}}(1+|t|_\infty^\beta)\right.\right.$$
$$\left.\left.\times\left(\left(\frac{1+|t|_\infty^{\lfloor d/2\rfloor+4}}{\sqrt{n}n^{1-\alpha}}\right)a_{n',p,\alpha,\beta'}(t')\right)^2 dt\right)^{1/2}\right)$$

$$= O\left(\frac{1}{\sqrt{n}n^{1-\alpha}}\sum_{l=0}^{\lfloor n/2\rfloor-\lceil\kappa n^\alpha\rceil-1}\left(\int_{|t|_\infty\le n^{1/2-\alpha}}(1+|t'|_\infty^{d+8+\beta})\right.\right.$$
$$\left.\left.\times(a_{n',p,\alpha,\beta'}(t'))^2 dt\right)^{1/2}\right)$$

$$= O\left(\frac{1}{\sqrt{n}n^{1-\alpha}}\sum_{l=0}^{\lfloor n/2\rfloor-\lceil\kappa n^\alpha\rceil-1}\left(\int_{|t'|_\infty\le n'^{1/2-\alpha}}(1+|t'|_\infty^{d+8+\beta})\right.\right.$$
$$\left.\left.\times(a_{n',p,\alpha,\beta'}(t'))^2 dt'\right)^{1/2}\right)$$

$$= O\left(\frac{1}{n^{2(1/2-\alpha)}}\right) = O\left(\frac{1}{n^{1/2-\alpha}}\right),$$

since $l\le\lfloor\frac{n}{2}\rfloor-\lceil\kappa n^\alpha\rceil-1$ implies $n'\ge\frac{n}{2}$.

(f) Using the fact that $(a_{m,p,\alpha,\beta'})_m$ is uniformly bounded, we have

$$\sum_{l=\lfloor n/2\rfloor-\lceil\kappa n^\alpha\rceil}^{n-\lfloor\kappa n^{\alpha/2}\rfloor-1}b_{n,l}(t)\left(\frac{1+|t|_\infty^3}{\sqrt{n}n^{1-\alpha}}\right)a_{n',p,\alpha,\beta'}(t')$$

$$= O\left(\sum_{l=\lfloor n/2 \rfloor - \lceil \kappa n^\alpha \rceil}^{n - \lfloor \kappa n^{\alpha/2} \rfloor - 1} \exp\left\{-\frac{\langle t, D(f)t \rangle}{2}\left(\frac{1}{2} - \frac{\kappa}{n^{1-\alpha}} - \frac{d+6}{2n}\right)\right\}\right.$$
$$\left. \times \left(\frac{1 + |t|_\infty^{\lfloor d/2 \rfloor + 4}}{\sqrt{n} n^{1-\alpha}}\right)\right)$$
$$= O\left(\exp\left\{-\frac{\langle t, D(f)t \rangle}{2}\left(\frac{1}{2} - \frac{\kappa}{n^{1-\alpha}} - \frac{d+6}{2n}\right)\right\}\left(\frac{1 + |t|_\infty^{\lfloor d/2 \rfloor + 4}}{n^{1/2-\alpha}}\right)\right).$$

Terms studied in (a) and (d)–(f) give contributions to the $a_{n,p+1,\alpha,\beta}$ term in (4) (for $p+1$ instead of $p$). Terms studied in (b) and (c) contribute to the first part of estimate (4) (for $p+1$ instead of $p$).

*Conclusion.* Now we deduce Theorem 2.2 from Proposition 2.7. Let a real number $\alpha \in ]0; \frac{1}{4}[$ and an integer $p \geq 2$ be given. Let us take $U_{n,p} := n^{(1/2-\alpha)(1-(1+d/2)/(p+d/2))}$. From Proposition 2.7, we get

$$\left(\int_{|t|_\infty \leq U_{n,p}} \sum_{k=0}^{\lfloor d/2 \rfloor + 1} \sum_{j_1,\ldots,j_k=1,\ldots,d} \left|\frac{\partial^k}{\partial t_{j_1} \cdots \partial t_{j_k}} h_n(f,t)\right|^2 dt\right)^{1/2}$$
$$= O_{n \to +\infty}\left(\frac{1}{n^{1/2-\alpha}}\right).$$

Finally, according to Yurinskii's result (recalled in Proposition 2.6 of this paper), we have

$$\forall \alpha \in \left]0; \frac{1}{4}\right[ \quad \forall p \geq 2 \quad \Pi_n(f) = O_{n \to +\infty}\left(\frac{1}{n^{(1/2-\alpha)(1-(1+d/2)/(p+d/2))}}\right).$$

**3. Limit theorem with rate of convergence for the averaging method.** We are interested in the asymptotic behavior of the error term between the solution of a differential equation perturbed by a transformation and the solution of the associated averaged differential equation. Results of convergence in distribution have been established in [20, 21, 28], for example.

3.1. *Averaging method for differential equation perturbed by a transformation.* In the following, we consider a (discrete-time) probability dynamical system $(\Omega, \mathcal{F}, \nu, T)$. Let an integer $d \geq 1$ be given. Let $F : \mathbf{R}^d \times \Omega \to \mathbf{R}^d$ be a measurable function uniformly bounded and uniformly Lipschitz continuous in the first parameter. We denote by $L_F$ its Lipschitz constant in the first parameter.

For any $\varepsilon > 0$ and any $(x, \omega)$ in $\mathbf{R}^d \times \Omega$, we consider the continuous solutions $(x_t^\varepsilon(x,\omega))_t$ and $(w_t(x))_t$ of the following differential equations (with

initial condition):

$$(24) \quad \forall\, t \in \mathbf{R} \setminus \varepsilon \mathbf{Z}, \qquad \frac{dx_t^\varepsilon}{dt}(x,\omega) = F(x_t^\varepsilon(x,\omega), T^{\lfloor t/\varepsilon \rfloor}(\omega)) \quad \text{and} \quad x_0^\varepsilon(x,\omega) = x$$

and

$$(25) \quad \frac{dw_t}{dt}(x) = \bar{F}(w_t(x)) = \int_\Omega F(w_t(x), \omega')\, d\nu(\omega') \quad \text{and} \quad w_0(x) = x.$$

Let us define the error term $(e_t^\varepsilon(x,\omega))_t$ as follows:

$$(26) \quad e_t^\varepsilon(x,\omega) := x_t^\varepsilon(x,\omega) - w_t(x).$$

NOTATION 3.1. Let a function $g: \mathbf{R}^d \times \Omega \to \mathbf{R}^d$ and an integer $k \geq 1$ be given.

We denote by $D_1^k g$ the $k$th differential of $g$ relative to the first parameter if it exists. Let us write $D_1 g := D_1^1 g$.

The function $g$ is said to be $C_b^{k,*}$ if $g$ is measurable, uniformly bounded, $C^k$-regular in the first parameter and if $D_1 g, \ldots, D_1^k g$ are measurable and uniformly bounded.

For any function $h: \mathbf{R}^d \to \mathbf{R}^d$, we denote by $D^k h$ the $k$th differential of $h$, if it is well defined. We write $Dh := D^1 h$.

We will make the following assumptions.

HYPOTHESIS 3.2. (i) The space $\Omega$ is endowed with a metric $d$, $\nu$ is a Borel measure ( for the topology induced by $d$ on $\Omega$) and there exists a real number $r_0 \geq 1$ such that the multiple decorrelation Property $(\mathcal{P}_{r_0})$ holds for $(\Omega, \mathcal{F}, \nu, T)$.

(ii) The function $F: \mathbf{R}^d \times \Omega \to \mathbf{R}^d$ is uniformly $\eta$-Hölder continuous in the second parameter.

(iii) The function $F: \mathbf{R}^d \times \Omega \to \mathbf{R}^d$ is $C_b^{2,*}$.

We will denote by $\widetilde{F}$ the function given by

$$\widetilde{F}(x,\omega) := F(x,\omega) - \bar{F}(x).$$

According to the proof of Theorem 2.1.3 of [28], we have the following result.

THEOREM 3.3. *Let a real number $T_0 > 0$ be given. Under Hypothesis 3.2, for any integer $L \geq 1$, we have*

$$\sup_{0<\varepsilon<1}\, \sup_{x \in \mathbf{R}^d}\, \sup_{0 \leq t \leq T_0} \left\| \frac{e_t^\varepsilon(x,\cdot)}{\sqrt{\varepsilon}} \right\|_L < +\infty.$$

MULTIPLE DECORRELATION AND CONVERGENCE RATE                29Moreover, for any $x \in \mathbf{R}^d$, the family of processes $((e_t^\varepsilon(x,\cdot))_{0 \le t \le T_0})_{\varepsilon > 0}$ converges in distribution [in $(\mathcal{C}([0,T_0]), \|\cdot\|_\infty)$ for measure $\nu$], when $\varepsilon$ goes to 0, to the Gaussian process $(e_t^0(x,\cdot))_{0 \le t \le T_0}$ solution of

$$e_t^0(x,\cdot) = v_t(x,\cdot) + \int_0^t D\bar{F}(w_s(x)) \cdot e_s^0(x,\cdot)\, ds,$$

where $v_t(x,\cdot)$ is a Gaussian process with independent increments, centered and such that

$$\mathbf{E}[(v_t(x,\cdot))^{\otimes 2}] = \int_0^t \mathcal{A}(\widetilde{F}(w_s(x),\cdot))\, ds,$$

with $\mathcal{A}(g) := \lim_{n \to +\infty} \mathbf{E}_\nu[(\frac{S_n(g)}{\sqrt{n}})^{\otimes 2}] = \mathbf{E}_\nu[g \otimes g] + \sum_{k \ge 1}(\mathbf{E}_\nu[g \otimes g \circ T^k] + \mathbf{E}_\nu[g \circ T^k \otimes g])$, for any $\nu$-centered, bounded $\eta$-Hlder continuous function $g : \Omega \to \mathbf{R}^d$.

An analogous result has been established in [21] under hypotheses of mixing for sub-$\sigma$-algebras (cf. also [20]).

### 3.2. Statement.

THEOREM 3.4. *Let $x \in \mathbf{R}^d$ and a real number $s > 0$ be given. Under Hypothesis 3.2, if $D_1 F$ is uniformly $\eta$-Hölder continuous in the second parameter, then the following limit exists:*

$$\Sigma_F^2 := \lim_{\varepsilon \to 0} \mathbf{E}_\nu\left[\left(\frac{e_s^\varepsilon(x,\cdot)}{\sqrt{\varepsilon}}\right)^{\otimes 2}\right].$$

*If, moreover, the matrixes $\mathcal{A}(\widetilde{F}(w_u(x),\cdot))$ defined above are nondegenerate (for all $u \in [0;s]$), then the family of random variables $(\frac{e_s^\varepsilon(x,\cdot)}{\sqrt{\varepsilon}})_{\varepsilon > 0}$ converges in distribution to a random variable with normal distribution $\mathcal{N}(0, \Sigma_F^2)$, and we have*

$$\forall \alpha > 0, \qquad \Pi\left(\nu_*\left(\frac{e_s^\varepsilon(x,\cdot)}{\sqrt{\varepsilon}}\right), \mathcal{N}(0, \Sigma_F^2)\right) = O(\varepsilon^{1/2 - \alpha}).$$

### 3.3. Proof.

Let us suppose $s = 1$ (this is not a restrictive hypothesis: it suffices to replace the function $F$ by the function $s \cdot F$). For any $(x,\omega) \in \mathbf{R}^d \times \Omega$ and any real number $\varepsilon > 0$, we define

$$v_t^\varepsilon(x,\omega) := \frac{1}{\sqrt{\varepsilon}} \int_0^t \widetilde{F}(w_s(x), T^{\lfloor s/\varepsilon \rfloor}(\omega))\, ds$$

and

$$y_t^\varepsilon(x,\omega) := \frac{1}{\sqrt{\varepsilon}} \int_0^t \exp\left\{\int_s^t D\bar{F}(w_r(x))\, dr\right\} \widetilde{F}(w_s(x), T^{\lfloor s/\varepsilon \rfloor}(\omega))\, ds.$$

$y_t^\varepsilon(x,\omega)$ is solution of $y_t^\varepsilon(x,\omega) = v_t^\varepsilon(x,\omega) + \int_0^t D\bar{F}(w_s(x)) \cdot y_s^\varepsilon(x,\omega)\,ds$. Our proof of Theorem 3.4 is based on the two following propositions (Propositions 3.5 and 3.7). The following result shows how the study of $\frac{e_1^\varepsilon(x,\cdot)}{\sqrt{\varepsilon}}$ comes down to the study of $y_1^\varepsilon(x,\cdot)$.

PROPOSITION 3.5. *Let a real number $T_0 > 0$ be given. Under Hypothesis 3.2, we have*

$$\forall p \in [1,+\infty[, \qquad \sup_{0 \le t \le T_0} \sup_{x \in \mathbf{R}^d} \left\| \frac{e_t^\varepsilon(x,\cdot)}{\sqrt{\varepsilon}} - y_t^\varepsilon(x,\cdot) \right\|_{L^p} = O(\varepsilon^{1/4}).$$

*If, moreover, function $D_1 F$ is uniformly $\eta$-Hölder continuous in the second variable, then we have*

$$\forall p \in [1,+\infty[, \qquad \sup_{0 \le t \le T_0} \sup_{x \in \mathbf{R}^d} \left\| \frac{e_t^\varepsilon(x,\cdot)}{\sqrt{\varepsilon}} - y_t^\varepsilon(x,\cdot) \right\|_{L^p} = O(\sqrt{\varepsilon}).$$

COROLLARY 3.6. *Under hypotheses of Theorem 3.4, we have*

$$\forall \alpha > 0, \qquad \mathcal{K}\left(\frac{e_t^\varepsilon(x,\cdot)}{\sqrt{\varepsilon}}, y_t^\varepsilon(x,\cdot)\right) = O(\varepsilon^{1/2-\alpha}).$$

PROOF. This is a consequence of the second point of Proposition 3.5. Effectively, if $X$ and $Y$ are two $\mathbf{R}^d$-valued random variables defined the a same probability space, then we have $\mathbf{P}(|X-Y|_\infty > \varepsilon) \le \frac{\|X-Y\|_p^p}{\varepsilon^p}$ and so $\mathcal{K}(X,Y) \le \|X-Y\|_p^{p/(p+1)}$.  □

PROOF OF PROPOSITION 3.5. The first point is a consequence of computations detailed in [28], Section 2.4, proof of Theorem 2.1.3 (cf. also [20], pages 220 and 221), these computations done in norm $L^1$ being still true in norm $L^p$ for any integer $p \ge 1$.

We only give the end of the proof of the second point which follows the scheme of the proof of the first point.

According to the computations done in [28], Section 2.4, identification of the cluster values, it is enough to show that we have

$$\sup_{0 \le t \le T_0} \sup_{x \in \mathbf{R}^d} \left\| \int_0^t D_1 \widetilde{F}(w_s(x), T^{\lfloor s/\varepsilon \rfloor}(\cdot)) \cdot y_s^\varepsilon(x,\cdot)\,ds \right\|_{L^p} = O(\sqrt{\varepsilon}),$$

for any integer $p \ge 1$. Let an integer $i = 1,\ldots,d$ be given. We have

$$\left(\int_0^t D_1 \widetilde{F}(w_s(x), T^{\lfloor s/\varepsilon \rfloor}(\cdot)) \cdot y_s^\varepsilon(x,\cdot)\,ds\right)_i = \sum_{j,k=1,\ldots,d} L_{i,j,k,\varepsilon}(t,x),$$

with

$$L_{i,j,k,\varepsilon}(t,x) = \varepsilon\sqrt{\varepsilon}\int_0^{t/\varepsilon}(D_1\widetilde{F}(w_{\varepsilon s}(x),T^{\lfloor s\rfloor}(\cdot)))_{i,j}$$

$$\times \left(\int_0^s\left(\exp\left\{\int_{\varepsilon u}^{\varepsilon s}D\bar{F}(w_r(x))\,dr\right\}\right)_{j,k}\right.$$

$$\left.\times \widetilde{F}_k(w_{\varepsilon u}(x),T^{\lfloor u\rfloor}(\cdot))\,du\right)ds.$$

Let $p$ be an even integer. We have

$$\|L_{i,j,k,\varepsilon}(t,x)\|_{L^p}^p$$

$$= \sqrt{\varepsilon}^p\varepsilon^p\mathbf{E}_\nu\left[\int_{B_{\varepsilon,p}}\left(\prod_{i'=1}^p(D_1\widetilde{F}(w_{\varepsilon s_{i'}}(x),T^{\lfloor s_{i'}\rfloor}(\cdot)))_{i,j}\right)\right.$$

$$\times\left(\prod_{j'=1}^p\left(\exp\left\{\int_{\varepsilon u_{j'}}^{\varepsilon s_{j'}}D\bar{F}(w_r(x))\,dr\right\}\right)_{j,k}\right.$$

$$\left.\left.\times \widetilde{F}_k(w_{\varepsilon u_{j'}}(x),T^{\lfloor u_{j'}\rfloor}(\cdot))\right)ds_1\cdots ds_p\,du_1\right.$$

$$\left.\cdots du_p\right]$$

$$= \sqrt{\varepsilon}^p\varepsilon^p\int_{B_{\varepsilon,p}}\left(\prod_{j'=1}^p\left(\exp\left\{\int_{\varepsilon u_{j'}}^{\varepsilon s_{j'}}D\bar{F}(w_r(x))\,dr\right\}\right)_{j,k}\right)$$

$$\times \mathbf{E}_\nu\left[\left(\prod_{i'=1}^p(D_1\widetilde{F}(w_{\varepsilon s_{i'}}(x),T^{\lfloor s_{i'}\rfloor}(\cdot)))_{i,j}\right)\right.$$

$$\left.\times\left(\prod_{j'=1}^p\widetilde{F}_k(w_{\varepsilon u_{j'}}(x),T^{\lfloor u_{j'}\rfloor}(\cdot))\right)\right]ds_1\cdots ds_p\,du_1\cdots du_p$$

$$\leq \sqrt{\varepsilon}^p\varepsilon^p e^{pT_0\|D\bar{F}\|_\infty}$$

$$\times \sum_{k_1,\ldots,k_{2p}=0}^{\lfloor T_0/\varepsilon\rfloor}\int_{k_{i'}\leq u_{i'}\leq k_{i'}+1}\left|\mathbf{E}_\nu\left[\prod_{i'=1}^{2p}G_{i'}(w_{\varepsilon u_{i'}}(x),T^{k_{i'}}(\cdot))\right]\right|du_1\cdots du_{2p},$$

with $B_{\varepsilon,p}:=\{(s_1,\ldots,s_p,u_1,\ldots,u_p)\in\mathbf{R}^{2p}:0\leq u_i\leq s_i\leq\frac{t}{\varepsilon}\}$ and by taking $G_{2i'-1}(x',\cdot)=(D_1\widetilde{F}(x',\cdot))_{i,j}$ and $G_{2i'}(x',\cdot)=\widetilde{F}_k(x',\cdot)$ for any $i'=1,\ldots,p$. According to Property $(\mathcal{P}_{r_0})$ and to the proof of Lemma 2.3.4 of [28], we know that, for any integer $L \geq 1$ and any real number $M > 0$, we have

$$\sup_{N \geq 1} \frac{1}{N^{L/2}} \sum_{n_1,\ldots,n_L=0}^{N-1} \sup_{H=(H^{(1)},\ldots,H^{(L)}) \in \mathcal{E}_{L,M}} \left| \mathbf{E}_\nu \left[ \prod_{k=1}^{L} H^{(i)} \circ T^{n_i} \right] \right| < +\infty,$$

where $\mathcal{E}_{L,M}$ is the set of $H = (H^{(1)}, \ldots, H^{(L)})$ where the functions $H^{(i)} : \Omega \to \mathbf{R}$ are bounded, $\eta$-Hölder continuous, $\nu$-centered and satisfy $\|H^{(i)}\|_\infty + C_{H^{(i)}}^{(\eta)} \leq M$. We get

$$\left\| \int_0^t D_1 \widetilde{F}(w_s(x), T^{\lfloor s/\varepsilon \rfloor}(\cdot)) y_s^\varepsilon(x, \cdot) \, ds \right\|_{L^p}^p = O(\varepsilon^{p/2}). \qquad \square$$

In the following, we study the behavior of the family of random variables $(y_1^\varepsilon(x, \cdot))_{\varepsilon > 0}$ when $\varepsilon$ goes to 0 (asymptotic behavior of the covariance matrices, convergence in distribution with rate of convergence). Let us notice that the study of the family of random variables $(y_1^\varepsilon(x, \cdot))_{\varepsilon > 0}$ when $\varepsilon$ goes to 0 comes down to the study of the sequence of random variables $(y_1^{1/N}(x, \cdot))_N$ when $N$ goes to $+\infty$. Effectively, we have

(27) $$\sup_{0 \leq s \leq T_0} \sup_{\omega \in \Omega} |y_s^\varepsilon(x, \omega) - y_s^{1/\lfloor 1/\varepsilon \rfloor}(x, \omega)|_\infty = O(\sqrt{\varepsilon}).$$

PROPOSITION 3.7.　*Under Hypothesis* 3.2, *the following limit exists*:

$$\Sigma_F^2 := \lim_{N \to +\infty} \mathbf{E}[(y_1^{1/N}(x, \cdot))^{\otimes 2}].$$

*If, moreover, the matrixes* $\mathcal{A}(\widetilde{F}(w(x), \cdot))$ *are nondegenerate* (*for all* $u \in [0; 1]$), *then we have*

$$\forall \alpha > 0, \quad \Pi(\nu_*(y_1^{1/N}(x, \cdot)), \mathcal{N}(0, \Sigma_F^2)) = O(N^{-1/2+\alpha}).$$

According to Proposition 3.5 and to (27), we have

$$\lim_{\varepsilon \to 0} \left\| \frac{e_t^\varepsilon(x, \cdot)}{\sqrt{\varepsilon}} - y_t^{1/\lfloor 1/\varepsilon \rfloor}(x, \cdot) \right\|_{L^2} = 0.$$

Hence, definitions of $\Sigma_F^2$ in Theorem 3.4 and in Proposition 3.7 coincide. Let us recall that, for any $\nu$-centered, bounded $\eta$-Hlder continuous function $g : \Omega \to \mathbf{R}^d$, we have defined

$$\mathcal{A}(g) := \mathbf{E}_\nu[g^{\otimes 2}] + \sum_{k \geq 1} (\mathbf{E}_\nu[g \otimes g \circ T^k] + \mathbf{E}_\nu[g \circ T^k \otimes g]).$$

LEMMA 3.8.　*Under Hypothesis* 3.2, *the following limit exists*:

$$\Sigma_F^2 := \lim_{N \to +\infty} \mathbf{E}_\nu[(y_1^{1/N}(x, \cdot))^{\otimes 2}]$$

*and satisfies*

$$\Sigma_F^2 := \frac{1}{N} \sum_{l=0}^{N-1} \mathcal{A}(F_{l,N}(w_{l/N}(x), \cdot)) + O\left(\frac{\log(N)^2}{N}\right),$$

*with*

$$F_{l,N}(x, \omega) := \int_0^1 \exp\left\{\frac{1}{N} \int_s^{N-l} D\bar{F}(w_{r/N}(x)) \, dr\right\} \widetilde{F}(w_{s/N}(x), \omega) \, ds.$$

PROOF. We have

$$y_1^{1/N}(x, \omega) = \frac{1}{\sqrt{N}} \int_0^N \exp\left\{\frac{1}{N} \int_s^N D\bar{F}(w_{r/N}(x)) \, dr\right\}$$

(28)
$$\times \widetilde{F}(w_{s/N}(x), T^{\lfloor s \rfloor}(\omega)) \, ds$$

$$= \frac{1}{\sqrt{N}} \sum_{k=0}^{N-1} F_{k,N}(w_{k/N}(x), T^k(\omega)).$$

Hence we have

$$\mathbf{E}_\nu[(y_1^{1/N}(x, \cdot))^{\otimes 2}] = \frac{1}{N} \sum_{k,l=0}^{N-1} \mathbf{E}_\nu[F_{k,N}(w_{k/N}(x), T^k(\cdot)) \otimes F_{l,N}(w_{l/N}(x), T^l(\cdot))].$$

We define

$$m_N := \frac{\log(N^{-2})}{\log(\delta_{r_0})} \quad \text{and} \quad A_N := \{(k, l) \in \{0, \ldots, N-1\}^2 : |k - l| \leq m_N\}.$$

We also define $B_N := \{0, \ldots, N-1\}^2 \setminus A_N$. According to the multiple decorrelation Property $(\mathcal{P}_{r_0})$ and to our choice of $m_N$, we have

$$\frac{1}{N} \sum_{(k,l) \in B_N} \mathbf{E}_\nu[F_{k,N}(w_{k/N}(x), T^k(\cdot)) \otimes F_{l,N}(w_{l/N}(x), T^l(\cdot))]$$

$$= O\left(\frac{1}{N} N^2 N^{-2}\right) = O\left(\frac{1}{N}\right).$$

On the other hand, since $\#A_N = O(Nm_N) = O(N\log(N))$, we have

$$\frac{1}{N} \sum_{(k,l) \in A_N} |\mathbf{E}_\nu[F_{k,N}(w_{k/N}(x), T^k(\cdot))$$

$$\otimes (F_{l,N}(w_{l/N}(x), T^l(\cdot)) - F_{k,N}(w_{k/N}(x), T^l(\cdot)))]|$$

$$= O\left(\frac{1}{N} N \log(N) \frac{\log(N)}{N}\right) = O\left(\frac{\log(N)^2}{N}\right).$$

Therefore, we have

$$\mathbf{E}_\nu[(y_1^{1/N}(x,\cdot))^{\otimes 2}]$$

$$= \frac{1}{N} \sum_{(k,l)\in A_N} \mathbf{E}_\nu[F_{k,N}(w_{k/N}(x), T^k(\cdot))$$

$$\otimes F_{k,N}(w_{k/N}(x), T^l(\cdot))] + O\left(\frac{\log(N)^2}{N}\right)$$

$$= \frac{1}{N} \sum_{k=m_N}^{N-1-m_N} \sum_{k-m_N \le l \le k+m_N} \mathbf{E}_\nu[F_{k,N}(w_{k/N}(x), T^k(\cdot))$$

$$\otimes F_{k,N}(w_{k/N}(x), T^l(\cdot))] + O\left(\frac{\log(N)^2}{N}\right)$$

$$= \frac{1}{N} \sum_{k=m_N}^{N-1-m_N} \mathcal{A}(F_{k,N}(w_{k/N}(x), \cdot)) + O\left(\frac{\log(N)^2}{N}\right)$$

$$= \frac{1}{N} \sum_{k=0}^{N-1} \mathcal{A}(F_{k,N}(w_{k/N}(x), \cdot)) + O\left(\frac{\log(N)^2}{N}\right). \qquad \square$$

PROOF OF PROPOSITION 3.7. The proof being analogous to the proof of Theorem 2.2 of the present paper, we do not give all its details. We only give the scheme of the $p$th iterative step. We will just detail computations which differ from the proof of Theorem 2.2. In the following, $N$ will be any integer and $t$ any point in $\mathbf{R}^d$ satisfying $|t|_\infty \le N^{1/2-\alpha}$.

Let us write $\Sigma_{F,l,N}^2 := \mathcal{A}(F_{l,N}(w_{l/N}(x), \cdot))$.

1. We define $H_0(t,N) := \exp\{\frac{-\langle t, \Sigma_F^2 t\rangle}{2}\} - \exp\{-\frac{1}{2N}\sum_{l=0}^{N-1}\langle t, \Sigma_{F,l,N}^2 t\rangle\}$. According to Lemma 3.8, there exists an integer $K_0 \ge 0$ such that we have

$$\sum_{k=0}^{\lfloor d/2 \rfloor+1} \sum_{j_1,\ldots,j_k=1,\ldots,d} \left|\frac{\partial^k}{\partial t_{j_1}\cdots\partial t_{j_k}} H_0(t,N)\right|$$
(29)
$$= O\left(\exp\left\{-\frac{1}{2}\left(\langle t, \Sigma_F^2 t\rangle - \frac{C^{te}(\log(N)^2)}{N}|t|_\infty^2\right)\right\}(1+|t|_\infty^{K_0})\frac{(\log(N))^2}{N}\right).$$

We define

$$H_1(t,N) := \exp\left\{-\frac{1}{2N}\sum_{l=0}^{N-1}\langle t, \Sigma_{F,l,N}^2 t\rangle\right\} - \prod_{l=0}^{N-1}\left(1 - \frac{1}{2N}\langle t, \Sigma_{F,l,N}^2 t\rangle\right).$$

We have

$$|H_1(t,N)| = O\left(\sum_{l'=0}^{N-1}\exp\left\{-\frac{1}{2N}\left\langle t, \left(\sum_{l\ne l'}\Sigma_{F,l,N}^2\right)t\right\rangle\right\}\frac{\sup_{l''}\langle t, \Sigma_{F,l'',N}^2 t\rangle^2}{8N^2}\right)$$

and, more generally, there exists a nonnegative integer $K_1$ such that

$$\sum_{k=0}^{\lfloor d/2\rfloor+1} \sum_{j_1,\ldots,j_k=1,\ldots,d} \left|\frac{\partial^k}{\partial t_{j_1}\cdots\partial t_{j_k}}H_1(t,N)\right|$$

$$(30) \qquad = O\left(\frac{1+|t|_\infty^{K_1}}{\sqrt{N}}\right.$$

$$\left.\times \exp\left\{-\frac{1}{2}\left(\left\langle t,\left(\sum_l \Sigma_{l,N}^2\right)t\right\rangle - \frac{d+4}{2N}\sup_{l''}\langle t,\Sigma_{l'',N}^2 t\rangle\right)\right\}\right).$$

We prove this estimate as we proved (6) in the proof of Theorem 2.2 by replacing (7) by the following formula which holds for any integer $N \geq 1$ and all $C^k$-regular functions $g_1,\ldots,g_N : \mathbf{R}^d \to \mathbf{C}$:

$$\frac{\partial^k}{\partial t_{j_1}\cdots\partial t_{j_k}}\left(\prod_{i=1}^N g_i\left(\frac{t}{\sqrt{N}}\right)\right)$$

$$= \sum_{m=1}^k \sum_{\{k_1,\ldots,k_m\}\in E_{m,N}} \sum_{\mathcal{A}\in\mathcal{L}_{m,k}} \left(\prod_{j\neq k_1,\ldots,k_m} g_j\left(\frac{t}{\sqrt{N}}\right)\right.$$

$$\left.\times \prod_{p=1}^m \frac{\partial^{\#A_p} g_{k_p}}{\partial t_{j_{l_1^{(p)}}}\cdots\partial t_{j_{l_{\#A_p}^{(p)}}}}\left(\frac{t}{\sqrt{N}}\right)\right)\frac{1}{N^{k/2}},$$

where $E_{m,N}$ is the set of subsets $\{1,\ldots,N\}$ with cardinal $m$ and where $\mathcal{L}_{m,k}$ is the set of partitions $\mathcal{A} = (A_1,\ldots,A_m)$ of $\{1,\ldots,k\}$ in nonempty subsets (i.e., $A_p \subseteq \{1,\ldots,k\}$, $A_p \neq \varnothing$, $\bigcup_p A_p = \{1,\ldots,k\}$ and $A_p \cap A_q = \varnothing$ if $p \neq q$) with $A_p = \{l_1^{(p)},\ldots,l_{\#A_p}^{(p)}\}$.

2. This leads us to the study of

$$\mathbf{E}_\nu\left[\exp\left\{\frac{i}{\sqrt{N}}\langle t, y_1^{1/N}(x,\cdot)\rangle\right\}\right] - \prod_{l=0}^{N-1}\left(1 - \frac{1}{2N}\langle t,\Sigma_{F,l,N}^2 t\rangle\right)$$

$$= \sum_{l=0}^{N-1}\left(\prod_{j=0}^{l-1}\left(1 - \frac{\langle t,\Sigma_{F,j,N}^2 t\rangle}{2N}\right)\right)$$

$$\times \mathbf{E}_\nu\left[Z_{l,N}(x,\cdot)\exp\left\{\frac{i}{\sqrt{N}}\sum_{k=l+1}^{N-1}\langle t,F_{k,N}(w_{k/N}(x),T^k(\cdot))\rangle\right\}\right],$$

with

$$Z_{l,N}(x,\cdot) := \exp\left\{\frac{i}{\sqrt{N}}\langle t,F_{l,N}(w_{l/N}(x),T^l(\cdot))\rangle\right\} - 1 + \frac{\langle t,\Sigma_{F,l,N}^2 t\rangle}{2N}.$$

3. We consider the quantities $M := p+3$ and $a_1, \ldots, a_M$ introduced in the proof of Theorem 2.2. We still define $A_0 := 0$ and $A_m := \sum_{j=1}^{m} a_j$ for every $m = 1, \ldots, M$. There exists a real number $\kappa > 0$ such that $a_j < \kappa N^{\alpha/2}$ for any $j = 1, \ldots, M$.

4. We estimate the following quantity as we have estimated $H_2$ in the proof of Theorem 2.2 [cf. estimate (10)]:

$$H_2(t,N) := \sum_{l=N-\lfloor \kappa N^{\alpha/2} \rfloor}^{N-1} \left( \prod_{q=0}^{l-1} \left( 1 - \frac{\langle t, \Sigma^2_{F,q,N} t \rangle}{2N} \right) \right)$$

$$\times \mathbf{E}_\nu \left[ Z_{l,N}(x,\cdot) \exp \left\{ \frac{i}{\sqrt{N}} \sum_{k=l+1}^{N-1} \langle t, F_{k,N}(w_{k/N}(x), T^k(\cdot)) \rangle \right\} \right].$$

5. For any $l \leq N - \lfloor \kappa N^{\alpha/2} \rfloor - 1$ and any $j = 1, \ldots, M$, we define

$$\mathcal{F}_j^{(l)} := \exp \left\{ \frac{i}{\sqrt{N}} \sum_{k=l+A_{j-1}+1}^{l+A_j} \langle t, F_{k,N}(w_{k/N}(x), T^k(\cdot)) \rangle \right\}$$

and

$$\mathcal{G}^{(l)} := \exp \left\{ \frac{i}{\sqrt{N}} \sum_{k=l+A_M+1}^{N-1} \langle t, F_{k,N}(w_{k/N}(x), T^k(\cdot)) \rangle \right\}.$$

We have

$$\mathbf{E}_\nu \left[ Z_{l,N}(x,\cdot) \exp \left\{ \frac{i}{\sqrt{N}} \sum_{k=l+1}^{N-1} \langle t, F_{k,N}(w_{k/N}(x), T^k(\cdot)) \rangle \right\} \right]$$

$$= \mathbf{E}_\nu \left[ Z_{l,N}(x,\cdot) \left( \prod_{j=1}^{M} \mathcal{F}_j^{(l)} \right) \mathcal{G}^{(l)} \right].$$

Moreover, as in the proof of Theorem 2.2, we can show that we have

$$\sum_{k=0}^{\lfloor d/2 \rfloor + 1} \sum_{j_1,\ldots,j_k=1,\ldots,d} \left| \frac{\partial^k}{\partial t_{j_1} \cdots \partial t_{j_k}} H_3(t,N) \right| = O\left( \frac{|t|_\infty^{p+1}}{N^{(1/2-\alpha)(p+1)}} \right),$$

with

$$H_3(t,N) := \sum_{l=0}^{N - \lfloor \kappa N^{\alpha/2} \rfloor - 1} \left( \prod_{j=0}^{l-1} \left( 1 - \frac{\langle t, \Sigma^2_{F,l,N} t \rangle}{2N} \right) \right)$$

$$\times \mathbf{E}_\nu \left[ Z_{l,N}(x,\cdot) \mathcal{F}_1^{(l)} \left( \prod_{j=2}^{M} (\mathcal{F}_j^{(l)} - 1) \right) \mathcal{G}^{(l)} \right].$$

6. It remains to estimate

$$\sum_{\varepsilon=(\varepsilon_1,\ldots,\varepsilon_M)} \mathbf{E}_\nu\left[Z_{l,N}(x,\cdot)\left(\prod_{i=1}^M \varepsilon_i\right)\mathcal{G}^{(l)}\right],$$

where the sum is taken over the $\varepsilon = (\varepsilon_1,\ldots,\varepsilon_M) \in \prod_{j=1}^M \{-1; \mathcal{F}_j^{(l)}\}$ with $\varepsilon_1 = \mathcal{F}_1^{(l)}$, the $\varepsilon_j$ being not all equal to $\mathcal{F}_j^{(l)}$. For any such vector $\varepsilon = (\varepsilon_1,\ldots,\varepsilon_{p+3})$, we define $j_0 := \max\{j \geq 2 : \varepsilon_j = -1\}$. We write

$$D_{l,\varepsilon}(N,t) := Z_{l,N}(x,\cdot)\prod_{j=1}^{j_0-1}\varepsilon_j$$

and

$$E_{l,\varepsilon}(N,t) := \left(\prod_{j=j_0+1}^M \mathcal{F}_j^{(l)}\right)\mathcal{G}^{(l)} = \exp\left\{\frac{it}{\sqrt{N}}\sum_{k=l+1+A_{j_0}}^{N-1} F_{k,N}(w_{k/N}(x), T^k(\cdot))\right\}.$$

In this study, we will use the following estimate instead of (11) (used in the proof of Theorem 2.2):

$$(31) \qquad \frac{\partial^k}{\partial t_{j_1}\cdots\partial t_{j_k}}\prod_{j=0}^{l-1}\left(1 - \frac{\langle t, \Sigma^2_{F,j,N}t\rangle}{2N}\right) = O(b_{N,l}(t)),$$

with

$$b_{N,l}(t) := \sum_{m=0}^{\min(\lfloor d/2\rfloor+1, l)} l\cdots(l-m+1)\frac{|t|_\infty^m}{N^m}$$

(32)

$$\times \exp\left\{-\frac{1}{2N}\left(\left\langle t, \sum_{j=0}^{l-1}\Sigma^2_{F,j,N}t\right\rangle - m\sup_{j''}\langle t, \Sigma^2_{F,j'',N}t\rangle\right)\right\}.$$

We will see that we have

$$\sum_{l=0}^{N-1} b_{N,l}(t) = O\left(\min\left(N, \frac{N}{|t|_\infty^2}\right)\right).$$

First, let us notice that there exists a real number $\tilde{c}_0 > 0$ such that, for all integers $N, L \geq 1$ and all $x \in \mathbf{R}^d$, we have

$$0 \leq \langle x, \Sigma^2_{F,L,N}x\rangle \leq \tilde{c}_0|x|_\infty^2.$$

On the other hand, since the symmetric matrices $\mathcal{A}(\widetilde{F}(w_n(x),\cdot))$ are nondegenerate, there exist an integer $N_1 \geq 1$ and a real number $\tilde{c}_1 > 0$ such that, for all integer $L \geq N_1$ and all $x \in \mathbf{R}^d$, we have

$$\left\langle x, \frac{1}{L}\sum_{l=0}^{L-1}\Sigma^2_{F,l,N}x\right\rangle \geq \tilde{c}_1|x|_\infty^2.$$

If $l \geq \max(N_1, \frac{2(\lfloor d/2 \rfloor + 1)\tilde{c}_0}{\tilde{c}_1})$, then we have

$$b_{N,l}(t) \leq \sum_{m=0}^{\min(\lfloor d/2 \rfloor + 1, l)} l \cdots (l-m+1) \frac{|t|_\infty^m}{N^m} \exp\left\{-\frac{l\tilde{c}_1}{4N}|t|_\infty^2\right\}.$$

Hence, we get

$$\sum_{l=\max(N_1, \lceil (2(\lfloor d/2 \rfloor + 1)\tilde{c}_0)/\tilde{c}_1 \rceil)}^{N-1} b_{N,l}(t) = O\left(\min\left(N, \frac{N}{|t|_\infty^2}\right)\right).$$

On the other hand, we have

$$\sum_{l \leq \max(N_1, \frac{2(\lfloor d/2 \rfloor + 1)\tilde{c}_0}{\tilde{c}_1})} b_{N,l}(t) = O(1).$$

*First step*: *estimate for the covariance.* We use Property $(\mathcal{P}_{r_0})$ as in the proof of Theorem 2.2 to estimate $\operatorname{Cov}_\nu(D_{l,\varepsilon}(N,t), E_{l,\varepsilon}(N,t))$.

*Second step*: *estimate for the first expectation.* We show that we have

$$(33) \quad \sum_{k=0}^{\lfloor d/2 \rfloor + 1} \sum_{j_1, \ldots, j_k = 1, \ldots, d} \left|\frac{\partial^k}{\partial t_{j_1} \cdots \partial t_{j_k}} \mathbf{E}_\nu[D_{l,\varepsilon}(N,t)]\right| = O\left(\frac{1 + |t|_\infty^3}{\sqrt{N} \cdot N^{1-\alpha}}\right).$$

Let us denote by $\mathcal{J}$ the following set:

$$\mathcal{J} := \{j = 1, \ldots, j_0 - 1 : \varepsilon_j = \mathcal{F}_j^{(l)}\}.$$

Let us recall that 1 belongs to $\mathcal{J}$. We have

$$|\mathbf{E}_\nu[D_{l,\varepsilon}(N,t)]|$$

$$= \left|\mathbf{E}_\nu\left[Z_{l,N}(x, \cdot) \exp\left\{\frac{i}{\sqrt{N}} \sum_{j \in \mathcal{J}} \sum_{k=l+A_{j-1}+1}^{l+A_j} \langle t, F_{k,N}(w_{k/N}(x), T^k(\cdot))t\rangle\right\}\right]\right|.$$

By noticing that we have

$$\sum_{j \in \mathcal{J}} \sum_{k=l+A_{j-1}+1}^{l+A_j} |(F_{k,N}(w_{k/N}(x), T^k(\cdot)) - F_{l,N}(w_{l/N}(x), T^k(\cdot)))|_\infty$$

$$= O\left(\frac{(\log(N))^2}{N}\right),$$

we are led to the study of $|\mathbf{E}_\nu[\widetilde{D}_{l,\varepsilon}(N,t)]|$, with

$$\widetilde{D}_{l,\varepsilon}(N,t) := Z_{l,N}(x, \cdot) \exp\left\{\frac{i}{\sqrt{N}} \sum_{j \in \mathcal{J}} \sum_{k=l+A_{j-1}+1}^{l+A_j} \langle t, F_{l,N}(w_{l/N}(x), T^k(\cdot))\rangle\right\}.$$

We can estimate this quantity as we have estimated the term $\mathbf{E}_\nu[D_{l,\varepsilon}(n,t)]$ appearing in the proof of Theorem 2.2. We will not rewrite all the computations. We will just detail the case $k=0$.

According to Taylor's formula, we get

$$Z_{l,N}(x,\cdot) = \frac{i}{\sqrt{N}}\langle t, F_{l,N}(w_{l/N}(x), T^l(\cdot))\rangle$$
$$+ \frac{1}{2N}\langle t, (\Sigma_{F,l,N}^2 - (F_{l,N}(w_{l/N}(x), T^l(\cdot)))^{\otimes 2})t\rangle + O\left(\frac{|t|_\infty^3}{N^{3/2}}\right)$$

and

$$\exp\left\{\frac{i}{\sqrt{N}}\sum_{j\in\mathcal{J}}\sum_{k=l+A_{j-1}+1}^{l+A_j}\langle t, F_{l,N}(w_{l/N}(x), T^k(\cdot))\rangle\right\}$$
$$= 1 + \frac{i}{\sqrt{N}}\sum_{j\in\mathcal{J}}\sum_{k=l+A_{j-1}+1}^{l+A_j}\langle t, F_{l,N}(w_{l/N}(x), T^k(\cdot))t\rangle + O\left(\frac{|t|_\infty^2}{N^{1-\alpha}}\right).$$

Therefore, we have

$$\mathbf{E}_\nu[\widetilde{D}_{l,\varepsilon}(N,t)]$$
$$= \frac{1}{2N}\mathbf{E}_\nu\left[\left\langle t, (\Sigma_{F,l,N}^2 - (F_{l,N}(w_{l/N}(x), T^l(\cdot)))^{\otimes 2}\right.\right.$$
$$- F_{l,N}(w_{l/N}(x), T^l(\cdot))$$
$$\left.\left.\otimes \sum_{j\in\mathcal{J}}\sum_{k=l+A_{j-1}+1}^{l+A_j} F_{l,N}(w_{l/N}(x), T^k(\cdot))\right)t\right\rangle\right]$$
$$+ O\left(\frac{|t|_\infty^3}{\sqrt{N}N^{1-\alpha}}\right)$$
$$= O\left(\frac{|t|_\infty^2}{N^2}\right) + O\left(\frac{|t|_\infty^3}{\sqrt{N}N^{1-\alpha}}\right) = O\left(\frac{1+|t|_\infty^3}{\sqrt{N}N^{1-\alpha}}\right).$$

*Third step: estimate for the second expectation.* We write $N' = N'_{N,l,\varepsilon} := N - (l+1) - A_{j_0}$ and $t' = t'_{N,l,\varepsilon} := t\sqrt{\frac{N'}{N}}$ et $\beta' := \beta + d + 8$. We estimate $\mathbf{E}_\nu[E_{l,\varepsilon}(N,t)]$ with the use of the inductive hypothesis as we have done in the proof of Theorem 2.2. Hence, we get

$$\sum_{k=0}^{\lfloor d/2\rfloor+1}\sum_{j_1,\ldots,j_k=1,\ldots,d}\left|\frac{\partial^k}{\partial t_{j_1}\cdots\partial t_{j_k}}\mathbf{E}_\nu[E_{l,\varepsilon}(N,t)]\right|$$

$$\leq O\left((1+|t|_\infty^{d/2+1})\exp\left\{-\frac{\tilde{c}_1|t|_\infty^2}{2}\left(1-\frac{l+1}{n}-\frac{\kappa}{n^{1-\alpha}}\right)\right\}\right)$$
$$+ L_{p,\alpha,\beta'}\frac{1+|t'|_\infty^p}{N'^{p(1/2-\alpha)}} + a_{N',p,\alpha,\beta'}(t').$$

Therefore, we got estimates analogous to those established in the proof of Theorem 2.2. We conclude in the same way with the use of (31) and (32). □

## APPENDIX

**Optimal and suboptimal estimates in norm $L^p$.** Let us consider a time-continuous dynamical system $(\mathcal{M},\mathcal{T},\mu,(Y_t)_{t\in\mathbf{R}})$, where $(\mathcal{M},\mathcal{T},\mu)$ is a probability space and where $(Y_t)_{t\in\mathbf{R}}$ is a family of $\mu$-preserving transformations of $M$ such that $(t,y)\mapsto Y_t(y)$ is measurable and satisfies $Y_0 = id$ and $Y_{t+s} = Y_t \circ Y_s$. Let us fix an integer $d \geq 1$.

Let us consider a measurable function $f:\mathbf{R}^d \times \mathcal{M} \to \mathbf{R}^d$ bounded, uniformly Lipschitz continuous in the first parameter such that, for any $(x,y) \in \mathbf{R}^d \times \mathcal{M}$, the functions $t \mapsto f(x,Y_t(y))$ are continuous on the right-hand side and limited on the left-hand side (i.e., they are cadlag functions), the set of discontinuity points being contained in a numerable set $\mathcal{D}_y$ independent of $x$. For all $\varepsilon > 0$ and all $(x,y) \in \mathbf{R}^d \times \mathcal{M}$, we consider the continuous piecewise $C^1$ function, $t \mapsto X_t^\varepsilon(x,y)$, solution of the following differential equation with initial condition:

$$(34)\quad X_0^\varepsilon(x,y)=x \quad\text{and}\quad \forall t\in\mathbf{R}\setminus\varepsilon\mathcal{D}_y, \quad \frac{dX_t^\varepsilon(x,y)}{dt}=f(X_t^\varepsilon(x,y),Y_{t/\varepsilon}(y)).$$

We are interested in the behavior of $(X_t^\varepsilon(x,y))_t$ when $\varepsilon$ goes to 0. We approximate $(X_t^\varepsilon(x,y))_t$ by the solution $(W_t(x))_t$ of the differential equation with initial condition obtained from (34) by averaging

$$(35)\quad W_0(x)=x \quad\text{and}\quad \forall t\in\mathbf{R}, \quad \frac{dW_t(x)}{dt}=\bar{f}(W_t(x)),$$

with $\bar{f}(x') := \int_\mathcal{M} f(x',y')\,d\mu(y')$.

This leads us to the study of the behavior of the error term $(E_t^\varepsilon(x,y))_t$ between the solution of the perturbed equation (34) and the solution of the equation (35) obtained by averaging

$$(36)\qquad\qquad E_t^\varepsilon(x,y) := X_t^\varepsilon(x,y) - W_t(x).$$

In [21] and [28], the question of convergence in distribution of $(\frac{E_t^\varepsilon(x,y)}{\sqrt{\varepsilon}})_t$ when $\varepsilon$ goes to 0 has been studied. The aim of this part is to establish estimates as optimal as possible of $\sup_{x\in\mathbf{R}^d}\|\sup_{0\leq t\leq T_0}|E_t^\varepsilon(x,\cdot)|_\infty\|_{L^p}$, with $p \in [1;+\infty]$.

In the following, we denote $\tilde{f}(x,y) := f(x,y) - \bar{f}(x)$.

If $\mathcal{M}$ is a compact manifold, if the flow $(Y_t)_t$ is $C^1$ and if $f$ is $C^1$ with compact support and satisfies the following condition of uniformly bounded variance:

$$(37) \qquad \sup_{x \in \mathbf{R}^d} \sup_{t > 0} \left\| \left| \frac{1}{\sqrt{t}} \int_0^t \tilde{f}(x, Y_s(\cdot)) \, ds \right|_\infty \right\|_{L^2} < +\infty,$$

Dumas and Golse established the following estimate (cf. [11]):

$$(38) \qquad \forall T_0 > 0, \qquad \int_{\mathbf{R}^d \times \mathcal{M}} \sup_{0 \leq t \leq T_0} |E_t^\varepsilon(x, \cdot)|_\infty \, dx \, d\mu(y) = O(\varepsilon^{1/3}).$$

Let us notice that their proof is still valid in the general context described at the beginning of this appendix, when $f : \mathbf{R}^d \times \mathcal{M} \to \mathbf{R}^d$ is a continuous function with compact support, $C_b^{1,*}$ (i.e., measurable, uniformly bounded, $C^1$ in the first variable with $D_1 f$ measurable and uniformly bounded) satisfying the following integrally bounded variance property:

$$(39) \qquad \int_{\mathbf{R}^d} \sup_{t > 0} \left\| \left| \frac{1}{\sqrt{t}} \int_0^t \tilde{f}(x, Y_s(\cdot)) \, ds \right|_\infty \right\|_{L^2} dx < +\infty$$

(cf. [27]). Let us notice that, $f$ having a compact support, condition (39) is weaker than condition (37).

In this appendix, we will make stronger hypotheses than conditions (39) and (37), which will enable us to establish estimates in $O(\varepsilon^{1/2})$ or in $O(|\ln(\varepsilon)|\varepsilon^{1/2})$ according to results due to Billingsley [5] and Serfling [35].

In Section A.1.1 we give optimal and suboptimal estimates for $\sup_{x \in \mathbf{R}^d} \| \sup_{0 \leq t \leq T_0} |e_t^\varepsilon(x, \cdot)|_\infty \|_{L^p}$ in the case of averaging method perturbed by a transformation (cf. Section 3.1). In Section A.1.4, we deduce from Section A.1.1 estimates for $\sup_{x \in \mathbf{R}^d} \| \sup_{0 \leq t \leq T_0} |E_t^\varepsilon(x, \cdot)|_\infty \|_{L^p}$ when the flow is associated (in some sense) to a transformation satisfying hypotheses of Section A.1.1.

**A.1.1. Perturbation by a transformation.** In the following, we are in the general context described at the beginning of Section 3.1 (before Hypothesis 3.2). We will suppose that this dynamical system is invertible, that is, that $T$ is one-to-one from a set $\Omega \setminus N_0$ onto a set $\Omega \setminus N_1$ with $\nu(N_0) = \nu(N_1) = 0$ and that the inverse transformation $T^{-1}$ is measurable. Such a hypothesis is not restrictive. Effectively, any dynamical system is a factor of an invertible dynamical system (its natural extension). We consider a real number $T_0 > 0$. We are interested in the study of the asymptotic behavior (as $\varepsilon$ goes to 0) of the following quantities:

$$(40) \qquad \sup_{x \in \mathbf{R}^d} \left\| \sup_{t \in [0; T_0]} |e_t^\varepsilon(x, \cdot)|_\infty \right\|_{L^p},$$

with $p \geq 1$. For any $(x, \omega) \in \mathbf{R}^d \times \Omega$, we define $\widetilde{F}(x, \omega) := F(x, \omega) - \bar{F}(x)$. According to Gronwall's lemma, we have

PROPOSITION A.1.1. *For any $\varepsilon > 0$ and any $(x,\omega) \in \mathbf{R}^d \times \Omega$, we have*

$$\sup_{t \in [0;T_0]} |e_t^\varepsilon(x,\omega)|_\infty \le (1 + L_F e^{L_F T_0}) \sup_{t \in [0;T_0]} \left| \varepsilon \int_0^{t/\varepsilon} \widetilde{F}(w_{\varepsilon s}(x), T^{\lfloor s \rfloor}(\omega)) \, ds \right|_\infty$$

*and*

$$\sup_{t \in [0;T_0]} \left| \varepsilon \int_0^{t/\varepsilon} \widetilde{F}(w_{\varepsilon s}(x), T^{\lfloor s \rfloor}(\omega)) \, ds \right|_\infty \le (1 + L_F T_0) \sup_{t \in [0;T_0]} |e_t^\varepsilon(x,\omega)|_\infty.$$

According to this result, the study of (40) brings us to the study of the following quantity:

$$\sup_{x \in \mathbf{R}^d} \left\| \sup_{t \in [0;T_0]} \left| \varepsilon \int_0^{t/\varepsilon} \widetilde{F}(w_{\varepsilon s}(x), T^{\lfloor s \rfloor}(\omega)) \, ds \right|_\infty \right\|_{L^p}.$$

A.1.2. *Estimate in norm $L^2$: a suboptimal result.* A first result is the following one.

THEOREM A.1.2. *If we have*

$$(41) \qquad \sup_{i=1,\ldots,d} \sum_{k \in \mathbf{Z}} \sup_{x,y \in \mathbf{R}^d} |\mathbf{E}_\nu[\widetilde{F}_i(x,\cdot) \cdot \widetilde{F}_i(y, T^k(\cdot))]| < +\infty,$$

*then we have*

$$(42) \qquad \sup_{x \in \mathbf{R}^d} \left\| \sup_{t \in [0;T_0]} |e_t^\varepsilon(x,\cdot)|_\infty \right\|_{L^2} = O(|\ln(\varepsilon)| \sqrt{\varepsilon}).$$

Let us notice that the condition (41) is close to the condition (37), the main difference being the fact that in (41) we study covariances of functions $\widetilde{F}_i(x,\cdot)$ and $\widetilde{F}_i(y,\cdot)$, with $x$ and $y$ maybe distinct. Condition (41) is not extremely restrictive; in particular, we can verify it for the examples studied in [11] without making more computations than those done to show that the condition (37) is satisfied.

Let us recall the following result.

THEOREM A.1.3 ([5], page 102). *Let two real numbers $\alpha \ge 1$ and $\beta \ge 1$ be given. Let $(X_n)_n$ be a sequence of real-valued random variables defined on the same probability space and a sequence of nonnegative real numbers $(u_n)_n$ such that, for all integer $n_0 \ge 0$ and $n \ge 1$, we have*

$$\mathbf{E}\left[ \left| \sum_{k=n_0}^{n_0+n-1} X_i \right|^\alpha \right] \le \left( \sum_{k=n_0}^{n_0+n-1} u_i \right)^\beta;$$

*then, for all integers $n_0 \geq 0$ and $n \geq 1$, we have*

$$\mathbf{E}\left[\sup_{m=1,\ldots,n}\left|\sum_{k=n_0}^{n_0+m-1} X_i\right|^\alpha\right] \leq (\log_2(4n))^\alpha \left(\sum_{k=n_0}^{n_0+n-1} u_i\right)^\beta.$$

SCHEME OF THE PROOF OF THEOREM A.1.2. Let us apply Theorem A.1.3 to $X_k := \int_k^{k+1} \widetilde{F}_i(w_{\varepsilon s}(x), T^{\lfloor s \rfloor}(\cdot)) \, ds$, $\alpha = 2$, $u_i = C$ and $\beta = 1$. We get

$$\sup_{x\in\mathbf{R}^d}\sup_{\varepsilon>0}\mathbf{E}_\nu\left[\sup_{n=0,\ldots,N}\left|\int_0^n \widetilde{F}(w_{\varepsilon s}(x), T^{\lfloor s \rfloor}(\cdot)) \, ds\right|_\infty^2\right] = O(N\log(N)).$$

We conclude with the use of the fact that $\widetilde{F}$ is uniformly bounded. □

The result of Theorem A.1.2 is suboptimal. Effectively, under hypotheses of Theorem A.1.2, we have

$$\sup_{x\in\mathbf{R}^d}\sup_{t\in[0;T_0]}\left\|\left|\int_0^{t/\varepsilon} \widetilde{F}(w_{\varepsilon s}(x), T^{\lfloor s \rfloor}(\cdot)) \, ds\right|_\infty\right\|_{L^2} = O(\varepsilon^{-1/2}).$$

If, moreover, we have $\sup_i \sum_{n\in\mathbf{Z}} |n| \sup_{x,y\in\mathbf{R}^d} |\mathbf{E}_\nu[\widetilde{F}_i(x,\cdot) \cdot \widetilde{F}_i(y, T^n(\cdot))]| < +\infty$, then a direct computation (cf. [20] and [28], Proposition 2.2.3) enables us to show that the covariance matrix (relative to $\nu$) of $\sqrt{\varepsilon}\int_0^{t/\varepsilon} \widetilde{F}(w_{\varepsilon s}(x), T^{\lfloor s \rfloor}(\cdot)) \, ds$ converges, as $\varepsilon$ goes to 0, to $\int_0^t \mathcal{A}(\widetilde{F}(w_u(x), \cdot)) \, du$, with $\mathcal{A}(g) = \sum_{k\in\mathbf{Z}} \mathbf{E}_\nu[g \otimes g \circ T^k]$. In that case, if some $\widetilde{F}(x,\cdot)$ are not coboundaries [i.e., if some matrices $\mathcal{A}(\widetilde{F}(w_u(x),\cdot))$ are not null], then $\sup_{x\in\mathbf{R}^d}\sup_{0\leq t\leq T_0}\mathbf{E}_\nu[|\int_0^{t/\varepsilon}\widetilde{F}(w_{\varepsilon s}(x), T^{\lfloor s \rfloor}(\cdot)) \, ds|_\infty^2]^{1/2}$ is exactly in $\frac{1}{\sqrt{\varepsilon}}$.

Therefore, according to Proposition A.1.1, $\sup_{x\in\mathbf{R}^d}\sup_{t\in[0;T_0]}\mathbf{E}_\nu[|e_t^\varepsilon(x,\cdot)|_\infty^2]^{1/2}$ is exactly in $\frac{1}{\sqrt{\varepsilon}}$.

Let us mention that the case when functions $\widetilde{F}(x,\cdot)$ are all coboundaries has been studied in [28].

Let us notice that we can get an estimate in $O(\sqrt{\varepsilon})$ in $L^2$ when we can apply the martingale method (see Gordin's method [15]; cf., e.g., Theorem 5.3.6 of [27]) with the use of Doob's inequality for martingales [16].

A.1.3. *Moment of larger order: optimal results.* We use the following result established in [35].

THEOREM A.1.4 (cf. Theorem B in [35]). *Let two real numbers $\alpha > 2$ and $C > 0$ be given. There exists a real number $K > 0$ such that, for any sequence of real random variables $(X_n)_n$ satisfying the following:*

$$\sup_{n_0\geq 0}\sup_{n\geq 1}\mathbf{E}_\nu\left[\frac{1}{n^{\alpha/2}}\left|\sum_{k=n_0}^{n_0+n-1} X_k\right|^\alpha\right] \leq C,$$

*we have*

$$\sup_{n_0 \geq 0} \sup_{n \geq 1} \mathbf{E}_\nu \left[ \frac{1}{n^{\alpha/2}} \sup_{k=1,\ldots,n} \left| \sum_{l=n_0}^{n_0+k-1} X_l \right|^\alpha \right] \leq K.$$

A consequence of this theorem is the following result.

THEOREM A.1.5. *Let an integer $p \geq 2$ be given. If the family of functions $\mathcal{F} := \{\widetilde{F}_i(x,\cdot); x \in \mathbf{R}^d, i = 1, \ldots, d\}$ satisfies the following condition:*

$$\sum_{l_1,\ldots,l_{2p}=0}^{N-1} \sup_{(g_1,\ldots,g_{2p}) \in \mathcal{F}^{2p}} \left| \mathbf{E}_\nu \left[ \prod_{i=1}^{2p} g_i \circ T^{l_i} \right] \right| = O(N^p),$$

*then we have*

$$\sup_{x \in \mathbf{R}^d} \left\| \sup_{t \in [0;T_0]} |e_t^\varepsilon(x,\cdot)|_\infty \right\|_{L^{2p}} = O(\sqrt{\varepsilon}).$$

PROOF. We have

$$\sup_{x \in \mathbf{R}^d} \sum_{l_1,\ldots,l_{2p}=0}^{N-1} \int_{l_1}^{l_1+1} \cdots \int_{l_{2p}}^{l_{2p}+1} \left| \mathbf{E}_\nu \left[ \prod_{j=1}^{2p} \widetilde{F}_i(w_{\varepsilon s_j}(x), T^{l_j}(\cdot)) \right] \right| ds_1 \cdots ds_{2p}$$
$$= O(N^p).$$

We conclude with Theorem A.1.4 for $X_k := \int_k^{k+1} \widetilde{F}_i(w_{\varepsilon s}(x), T^{\lfloor s \rfloor}(\cdot)) \, ds$ and for $\alpha = 2p$. □

Examples of systems satisfying the hypotheses of Theorem A.1.5 for all $p \in [1, +\infty[$ have been studied in [28]. In particular, we have the following result.

PROPOSITION A.1.6. *Under the two first points of Hypothesis 3.2, we have*

$$\forall p \in [1, +\infty[ \qquad \sup_{x \in \mathbf{R}^d} \left\| \sup_{t \in [0;T_0]} |e_t^\varepsilon(x,\cdot)|_\infty \right\|_{L^{2p}} = O(\sqrt{\varepsilon}).$$

PROOF. By a combinatorial argument (cf., e.g., the proof of Lemma 2.3.4 of [28]), we can show that, in this situation, hypotheses of Theorem A.1.5 are satisfied. □

**A.1.4. Perturbation by a flow.** We study here quantities $\sup_{x \in \mathbf{R}^d} \| \sup_{t \in [0;T_0]} |E_t^\varepsilon(x,\cdot)|_\infty \|_{L^p}$ for the averaging method for differential equations perturbed by a flow in the context described at the beginning of the Appendix. We will see how we can be brought to the question of the study of $\sup_{x \in \mathbf{R}^d} \| \sup_{t \in [0;T_0]} |e_t^\varepsilon(x,\cdot)|_\infty \|_{L^p}$, where $e_t^\varepsilon(x,\omega)$ is the error term in the averaging method for a differential equation perturbed by a transformation. We will consider the transformation $T = Y_1$ (in the case of diagonal flows) or we will use a representation of the flow as a special flow (in the case of the billiard flow). We will conclude with the help of the results of Section A.1.1.

A.1.5. *Flow stopped at time* 1. In this section, we take $(\Omega, \mathcal{F}, \nu, T) = (\mathcal{M}, \mathcal{T}, \mu, Y_1)$. We consider the function $F: \mathbf{R}^d \times \Omega \to \mathbf{R}^d$ defined by $F(x,\omega) := \int_0^1 f(x, Y_s(\omega))\, ds$. We consider the processes $(x_t^\varepsilon(x,\omega))$, $(w_t(x))$ and $(e_t^\varepsilon(x,\omega))$ given by (24), (25) and (26) for this choice of $(\Omega, \mathcal{F}, \nu, T)$ and of $F$. Then we can show that, for any real number $T_0 > 0$, we have

$$\sup_{(x,\omega) \in \mathbf{R}^d \times \Omega} \sup_{t \in [0,T_0]} |E_t^\varepsilon(x,\omega) - e_t^\varepsilon(x,\omega)|_\infty = O(\varepsilon).$$

According to results established in [23] about diagonal flows, Property $(\mathcal{P}_1)$ is satisfied in this context. This enables us to show the following result, according to Proposition A.1.6.

EXAMPLE A.1.7 (Diagonal flow on a homogeneous space). Let $d \geq 2$ be an integer and let $\Gamma$ be a cocompact subgroup of $G := SL(d, \mathbf{R})$. We consider the quotient space $\mathcal{M} := SL(d, \mathbf{R})/\Gamma$ endowed with the probability measure left-translation-invariant $\bar{\mu}$ induced on $\mathcal{M}$ by the Haar measure on $G$. Let $(T_i)_{i=1}^d$ be a decreasing sequence of $d$ positive real numbers not all equal to 1, the product of which is 1. For any real number $t \in \mathbf{R}$, we denote by $T^t$ the matrix

$$T^t = \begin{pmatrix} T_1^t & & & & \\ & T_2^t & & 0 & \\ & & \ddots & & \\ & 0 & & T_{d-1}^t & \\ & & & & T_d^t \end{pmatrix}.$$

We consider the $\bar{\mu}$-preserving flow $(Y_t)_t$ defined on $G/\Gamma$ by $Y_t(x\Gamma) = T^t x \Gamma$. Let us fix a riemannian metric $d_0$ on $G$ invariant by right-translation and let us define a metric $d$ on $G/\Gamma$ by

$$d(x\Gamma, y\Gamma) := \inf_{\gamma \in \Gamma} d_0(x, y\gamma).$$

If $f: \mathbf{R}^d \times \mathcal{M} \to \mathbf{R}^d$ is a measurable function, which is uniformly bounded and uniformly Lipschitz continuous in the first variable and uniformly Hlder in the second variable (for the metric $d$), then we have

$$\forall p \in [1, +\infty[, \qquad \sup_{x \in \mathbf{R}^d} \left\| \sup_{t \in [0;T_0]} |E_t^\varepsilon(x, \cdot)|_\infty \right\|_{L^{2p}} = O(\sqrt{\varepsilon}),$$

where $(X_t^\varepsilon(x, \omega))$, $(W_t(x))$ and $(E_t^\varepsilon(x, \omega))$ are defined by (34), (35) and (36) for this choice of $(\mathcal{M}, \mathcal{T}, \mu, (Y_t)_t)$ and of $f$.

A.1.6. *Case of a special flow.* Let us suppose now that the time-continuous dynamical system $(\mathcal{M}, \mathcal{T}, \mu, (Y_t)_t)$ is the special flow associated to a dynamical system $(\Omega, \nu, T)$ and to a roof function $\tau: \Omega \to [0; +\infty[$ satisfying $\inf \tau > 0$ and $\sup \tau < +\infty$, which means:

(i) $\mathcal{M}$ is the set $\{(\omega, s): \omega \in \Omega, s \in [0; \tau(\omega)]\}$ with identifications $(\omega, \tau(\omega)) \equiv (T(\omega), 0)$;

(ii) $\mathcal{T}$ is the $\sigma$-algebra induced on $\mathcal{M}$ by the product $\sigma$-algebra $\Omega \times \mathbf{R}_+$;

(iii) the probability measure $\mu$ is given by $d\mu(\omega, s) \frac{1}{\int_\Omega \tau \, d\nu} \, d\nu(\omega) \, ds$;

(iv) the flow $(Y_t)_t$ is given by $Y_t(\omega, s) = (\omega, s+t)$ with the identifications $(\omega, \tau(\omega)) \equiv (T(\omega), 0)$.

We make the following hypothesis on the function $f: \mathbf{R}^d \times \mathcal{M} \to \mathbf{R}^d$:

HYPOTHESIS A.1.8. The function $f$ is measurable, uniformly bounded and uniformly Lipschitz continuous in the first variable.

For every $(x, \omega) \in \mathbf{R}^d \times \Omega$, the function $s \mapsto f(x, (\omega, s))$ is continuous on $[0, \tau(\omega)[$ and the following limit exists: $\lim_{s \to \tau(\omega)-} f(x, (\omega, s))$.

We then consider the function $F: \mathbf{R}^d \times \Omega \to \mathbf{R}^d$ defined by $F(x, \omega) := \int_0^{\tau(\omega)} f(x, (\omega, s)) \, ds$. We consider also the processes $(x_t^\varepsilon(x, \omega))$, $(w_t(x))$ and $(e_t^\varepsilon(x, \omega))$ defined by the (24), (25) and (26) for $(\Omega, \mathcal{F}, \nu, T)$ and for this choice of $F$. We consider the process $(f_t^\varepsilon(x, \omega))$ defined as the process $(e_t^\varepsilon(x, \omega))$ by replacing $F$ by the function $G$ given by $G(x, \omega) = \tau(\omega) \bar{f}(x)$. According to [21] (see also [28], Section 3.2), we have:

REMARK A.1.9. Under Hypothesis A.1.8, for any real number $T_0 > 0$, we have:

$$\sup_{x \in \mathbf{R}^d} \sup_{\varepsilon > 0} \sup_{t \in [0, T_0]} |E_t^\varepsilon(x, \omega) - (e_{\varepsilon n(t/\varepsilon, \omega)}^\varepsilon(x, \omega) - f_{\varepsilon n(t/\varepsilon, \omega)}^\varepsilon(x, \omega))|_\infty = O(\varepsilon),$$

with $n(t, \omega) := \max\{n \geq 0: \sum_{k=0}^{n-1} \tau(T^k(\omega)) \leq t\}$.

Hence, for any $T_0 > 0$ and any $p \in [1; +\infty[$, we have:

$$\left\| \sup_{t \in [0;T_0]} |E_t^\varepsilon(x, \cdot)|_\infty \right\|_{L^p} \leq \left\| \sup_{t \in [0;T_0/\inf_\Omega \tau]} |e_t^\varepsilon(x, \cdot) - f_t^\varepsilon(x, \cdot)|_\infty \right\|_{L^p} + O(\varepsilon).$$

On the other hand, as for Proposition A.1.1, we can show that we have:

REMARK A.1.10. Under Hypothesis A.1.8, there exists a real number $C > 0$ such that, for any real number $\varepsilon > 0$ and any $(x, \omega) \in \mathbf{R}^d \times \Omega$, we have

$$\sup_{t \in [0;T_0]} |e_t^\varepsilon(x,\omega) - f_t^\varepsilon(x,\omega)|_\infty \leq C \sup_{t \in [0;T_0]} \left| \varepsilon \int_0^{t/\varepsilon} H(w_{\varepsilon s}(x), T^{\lfloor s \rfloor}(\omega)) \, ds \right|_\infty,$$

with $H(x,\omega) := F(x,\omega) - \tau(\omega)\bar{f}(x)$.

According to the results on the billiard flow established in [28] (cf. also [37] and [38]), Property $(\mathcal{P}_r)$ is satisfied for every real number $r > 1$. Therefore, according to the proof of Theorem A.1.5, we have:

EXAMPLE A.1.11 (Billiard flow with finite horizon). Let $Q$ be a compact subset of the torus $\mathbf{T}^2 = \frac{\mathbf{R}^2}{\mathbf{Z}^2}$, the complement of which is a finite union of strictly convex open sets (open disks, e.g.) with closure pairwise disjoint and the boundary of which is $C^3$ with curvature never null. We are interested in the behavior of a point particle moving in $Q$ with unitary speed and elastic reflections off $\partial Q$. We consider the time-continuous dynamical system $(\mathcal{M}, \mathcal{T}, \mu, (Y_t)_t)$ defined as follows:

(a) we denote by $\mathcal{T}^1 Q$ the set of position-speed couples $(q, \vec{v})$ with $q \in Q$ and $\|\vec{v}\| = 1$; we define $\mathcal{M} := \{(q, \vec{v}) \in \mathcal{T}^1 Q : q \notin \partial Q \text{ or } \langle \vec{n}(q), \vec{v} \rangle \geq 0\}$, where $\vec{n}(q)$ is the unitary normal vector to $\partial Q$ in $q$ (oriented to the inside of $Q$) if $q \in \partial Q$. We endow $\mathcal{M}$ with the metric $d$ given by

$$d((q, \vec{v}), (q', \vec{v}')) = d_0(q, q') + d_1(\vec{v}, \vec{v}'),$$

where $d_0$ is the metric induced on $\mathbf{T}^2$ by the usual euclidean metric on $\mathbf{R}^2$ and where $d_1(\vec{v}, \vec{v}')$ is the absolute value of the angular measure taken in $]-\pi; \pi]$ of the angle $\widehat{(\vec{v}, \vec{v}')}$;

(b) $\mu$ is the normalized Lebesgue measure on $\mathcal{M}$;

(c) $(Y_t)_t$ is the billiard flow defined on $\mathcal{M}$ by $Y_t(q, \vec{v}) = (q', \vec{v}')$ is the position-speed couple at time $t$ of a particle that was at position $q$ with speed $\vec{v}$ at time 0.

For every $(q, \vec{v}) \in \mathcal{M}$, we define $\tau^+(q, \vec{v}) := \inf\{s > 0 : q + s\vec{v} \in \partial Q\}$. Let us suppose that function $\tau^+$ is bounded (we say that the billiard has finite horizon). If $f : \mathbf{R}^d \times \mathcal{M} \to \mathbf{R}^d$ is a measurable function, uniformly bounded, uniformly Lipschitz in the first variable and uniformly Hlder in the second variable, then we have

$$\forall p \in [1, +\infty[, \quad \sup_{x \in \mathbf{R}^d} \left\| \sup_{t \in [0;T_0]} |E_t^\varepsilon(x, \cdot)|_\infty \right\|_{L^{2p}} = O(\sqrt{\varepsilon}),$$

where the processes $(X_t^\varepsilon(x,\omega))$, $(W_t(x))$ and $(E_t^\varepsilon(x,\omega))$ have been defined by (34), (35) and (36) for this choice of $(\mathcal{M}, \mathcal{T}, \mu, (Y_t)_t)$ and of $f$.

**Acknowledgments.** I am grateful to J.-P. Conze for all the interesting discussions we have had about the subject of this article. I would also like to thank B. Courbot for his explanations about the speed of convergence in probabilistic limit theorems and for indicating to me the work of Yurinskii [39].

DEPARTMENT DE MATHÉMATIQUES
UNIVERSITÉ DE BREST
29285 BREST CEDEX
FRANCE
E-MAIL: [francoise.pene@univ-brest.fr](francoise.pene@univ-brest.fr)